\documentclass[12pt,doublespace]{article}

\newcommand{\pr}{{\it Proof.} \ }
\newcommand{\cc}{\frac{1}{2}}

\usepackage{graphicx}
\usepackage{amsfonts}
\usepackage{amsmath}
\usepackage{amssymb}
\usepackage{latexsym}

\begin{document}
\title{The distribution of the summatory function of the M\"{o}bius function}
\author{Nathan Ng}
\maketitle

\begin{abstract}
Let the summatory function of the M\"{o}bius function be denoted
$M(x)$. We deduce in this article conditional results concerning
$M(x)$ assuming the Riemann Hypothesis and a conjecture of Gonek
and Hejhal on the negative moments of the Riemann zeta function.
The main results shown are that the weak Mertens conjecture and
the existence of a limiting distribution of $e^{-y/2}M(e^{y})$ are
consequences of the aforementioned conjectures.  By probabilistic
techniques, we present an argument that suggests $M(x)$ grows as
large positive and large negative as a constant times $\pm
\sqrt{x} (\log \log \log x)^{\frac{5}{4}}$ infinitely often, thus
providing evidence for an unpublished conjecture of Gonek's.
\end{abstract}

\section{Introduction}

\footnotetext{Mathematics Subject Classification 11M26,11N56}
The M\"{o}bius function is defined
for positive integers $n$ by
\begin{equation}
 \mu(n) = \left\{ \begin{array}{ll}
                  1 & \mbox{if $n=1$} \\
                  0 & \mbox{if $n$ is not squarefree} \\
                  (-1)^{k} & \mbox{if $n$ is squarefree
                  and $n= p_{1} \ldots p_{k}$}
                  \end{array}
          \right.   \ .
\end{equation}
Its summatory function $M(x) = \sum_{n \le x} \mu(n)$ is closely related
to the reciprocal of the Riemann zeta function.
This connection may be observed by the identities
\[
   \frac{1}{\zeta(s)} = \sum_{n = 1}^{\infty} \frac{\mu(n)}{n^{s}}
   = s \int_{1}^{\infty} \frac{M(x)}{x^{s+1}} \ dx
     \label{eq:i1}
\]
valid for $\mathrm{Re}(s) > 1$ and
\begin{equation}
   M(x) = \frac{1}{2 \pi i}
   \int_{c - i \infty}^{c + i \infty}
   \frac{x^{s}}{s \zeta(s)} \ ds
   \label{eq:i2}
\end{equation}
where $c > 1$ and $x \not \in \mathbb{Z}$. In the theorems of this
article, we assume the truth of the Riemann hypothesis (RH) which
asserts that all non-real zeros of $\zeta(s)$ take the form $\rho
= \frac{1}{2} + i \gamma$ with $\gamma \in \mathbb{R}$. At times,
we also assume that all zeros of the zeta function are simple. It
is widely expected that all zeros of the zeta function are simple.
Currently, the best unconditional result is that at least
two-fifths of the zeros are simple \cite{C}. In light
of~(\ref{eq:i2}), sums of the form
\[ J_{-k}(T) =
   \sum_{0 < \gamma \le T} \frac{1}{|\zeta^{'}(\rho)|^{2k}}
\]
where $k \in \mathbb{R}$ are important in obtaining information
concerning $M(x)$. From different points of view Gonek \cite{G}
and Hejhal \cite{Hej} independently conjectured that
\begin{equation}
   J_{-k}(T) \asymp T (\log T)^{(k-1)^{2}} .
   \label{eq:gh}
\end{equation}
Gonek studied Dirichlet polynomial approximations of these
moments, whereas Hejhal studied the value distribution of $\log
\zeta^{'}(\rho)$ employing ideas of Selberg's. Henceforth, the
former assumption~(\ref{eq:gh}) will be referred to as the
Gonek-Hejhal conjecture. For $k=0$ we have $J_{0}(T)=N(T)$ where
$N(T)$ is the number of zeros in the box with vertices $0,1,1+iT$,
and $iT$.  Von Mangoldt (see \cite{D} pp. 97-100) proved that
\begin{equation}
   J_{0}(T) = \frac{T}{2 \pi} \log \frac{T}{2 \pi e}
   + O( \log T) \ .
   \label{eq:nt}
\end{equation}
For $k=1$ Gonek \cite{G2} conjectured the asymptotic formula
\begin{equation}
  J_{-1}(T) \sim \frac{3}{\pi^{3}}T \ .
  \label{eq:gm1}
\end{equation}
Moreover, he proved that $J_{-1}(T) \gg T$ (see \cite{G}) subject
to RH and all zeros of the Riemann zeta function are simple.
Recently, Hughes et al.\ \cite{HKO} using random matrix model
techniques conjectured that
\begin{equation}
   \sum_{0 < \gamma \le T} | \zeta^{'}(\rho) |^{2k}
   \sim
   \frac{G^{2}(k+2)}{G(2k+3)}a_{k} \frac{T}{2 \pi}
   \left( \log \frac{T}{2 \pi} \right)^{(k+1)^{2}}
   \label{eq:rmc}
\end{equation}
for $k > - \frac{3}{2}$ where
\[ a_{k} = \prod_{p} \left( 1 - \frac{1}{p} \right)^{k^{2}}
   \left( \sum_{m=0}^{\infty} \left( \frac{\Gamma(m+k)}
   {m! \Gamma(k)} \right)^{2} p^{-m} \right)
\]
and $G$ is Barnes' function defined by
\[ G(z+1) = (2 \pi)^{z/2} \exp \left( -\frac{1}{2}(z^{2} +
   \gamma z^{2} + z ) \right) \prod_{n=1}^{\infty}
    \left( \left(1 + \frac{z}{n} \right)^{n}
   e^{-z + z^{2}/2n} \right)
\]
where $\gamma$ denotes Euler's constant. One should note that in
the above definition of $a_{k}$, one must take an appropriate
limit if $k=0$ or $k=-1$.  Furthermore, one may check that $G(1) =
1$ and $a_{-1} = \frac{6}{\pi^{2}}$ and hence~(\ref{eq:rmc})
implies~(\ref{eq:gm1}) and moreover
it agrees with~(\ref{eq:nt}). \\
One notes that Gonek \cite{G2} arrives at
conjecture~(\ref{eq:gm1}) by pursuing ideas of Montgomery's
concerning the zero spacings (pair-correlation) of the zeta
function. On the other hand, the random matrix technique
originated with the work of Keating and Snaith \cite{KS}.  Their
idea was to model the Riemann zeta function by the characteristic
polynomial of a large random unitary matrix. They computed moments
of these characteristic polynomials averaged over the group of
unitary matrices.  These moments are much simpler to evaluate
since they may be transformed into the well-studied Selberg
integral. Following the work of Keating and Snaith, other authors
have used this analogy to speculate on the exact nature of certain
families of $L$-functions.  This analogy has been viewed as rather
fruitful, since to date it has always produced conjectures that
agree with known theorems.

In this article we deduce results about $M(x)$ assuming the
Riemann Hypothesis and the conjectural bound
\begin{equation}
   J_{-1}(T) = \sum_{0 < \gamma < T} \frac{1}{|\zeta^{'}(\rho)|^{2}}
   \ll T \ .
   \label{eq:j1}
\end{equation}
By making assumption~(\ref{eq:j1}), we implicitly assume that all
zeros are simple. If there were a multiple zero of $\zeta(s)$,
$J_{-1}(T)$ would be undefined for sufficiently large $T$
and~(\ref{eq:j1}) would fail to make sense. For a long time,
number theorists were interested in $M(x)$ as RH was a consequence
of the famous Mertens conjecture which states that
\[
   |M(x)| \le x^{\frac{1}{2}} \ \mathrm{for} \ x \ge 1 .
\]
For an excellent historical account of work on this problem see \cite{OR}.
An averaged version of this conjecture is the weak Mertens conjecture
which asserts that
\begin{equation}
   \int_{2}^{X} \left( \frac{M(x)}{x} \right)^{2} \, dx
   \ll \log X \ .
    \label{eq:wmc}
\end{equation}
The weak Mertens conjecture implies RH, all zeros of $\zeta(s)$
are simple, and that $\sum_{\gamma > 0} \frac{1}{|\rho
\zeta^{'}(\rho)|^{2}}$ converges. These consequences are proven in
Titchmarsh \cite{T} pp.\,376-380. Not surprisingly, the Mertens
conjecture was disproven by Odlyzko and te Riele \cite{OR} as they
showed that
\[ \liminf_{x \to \infty}  \frac{M(x)}{\sqrt{x}} < -1.009 \
   \mathrm{and} \
   \limsup_{x \to \infty} \frac{M(x)}{\sqrt{x}} > 1.06 \ .
\]
However, they did not actually provide a specific counterexample
to~(\ref{eq:mc}). In fact, the Mertens conjecture was put in
serious doubt many years earlier when Ingham \cite{In} proved
\[ \liminf_{x \to \infty} \frac{M(x)}{\sqrt{x}} = - \infty \
   \mathrm{and} \
   \limsup_{x \to \infty}  \frac{M(x)}{\sqrt{x}} = \infty
\]
assuming the following conjecture:

\noindent {\bf Linear independence conjecture} (LI) Assume
$\zeta(s)$ satisfies the Riemann Hypothesis.  If its zeros are
written as $\frac{1}{2} + i \gamma$, then the positive imaginary
ordinates of the zeros are linearly independent over $\mathbb{Q}$.

Currently there is very little numerical or theoretical evidence
supporting this conjecture. However, it is considered rather
unlikely that the imaginary ordinates of the zeros of the zeta
function satisfy any linear relations. The linear independence
conjecture has been used in the past to get a handle on some very
difficult problems in number theory (see \cite{In}, \cite{Mo},
\cite{RS}). For some modest numerical computations see \cite{B}.
Despite the above results, we would like to have a better
understanding of what the upper and lower bounds of
$x^{-\frac{1}{2}} M(x)$ should be. The true order of $M(x)$ is
something of a mystery. In fact, Odlyzko and te Riele \cite{OR}
p.\,3 comment that ``No good conjectures about the rate of growth
of $M(x)$ are known.'' Motivated by this comment, we attempt to
give an explanation of the true behaviour of $M(x)$ assuming
reasonable conjectures about the zeta function.

We briefly mention some notation used  throughout this article.
We will denote a sequence of effectively computable
positive constants as $c_{1},c_{2},c_{3}, \ldots$.
We will also employ the following notation.
Let $f(x),g(x)$ be two real valued functions with $g(x) >0$.
Then the notation $f(x) = \Omega_{+}(g(x))$ means
\[
   \limsup_{x \to \infty} \frac{f(x)}{g(x)} > 0
\]
and $f(x) = \Omega_{-}(g(x))$ means
\[
   \liminf_{x \to \infty} \frac{f(x)}{g(x)} < 0 \ .
\]

We now state our current knowledge of $M(x)$.
The best known unconditional upper bound is
\[ M(x) = O \left( x \exp \left( -c_{1} \log^{\frac{3}{5}} x
            (\log \log x)^{-\frac{1}{5}} \right) \right)
\]
for $c_{1} >0$ (see Ivi\'{c} \cite{Ivic} pp.\ 309-315) .
However, the Riemann hypothesis is equivalent to the bound
\[ M(x) = O\left( x^{\frac{1}{2}}
   \exp \left( \frac{c_{2} \log x}
   {\log \log x} \right) \right) \]
for $c_{2}>0$ (see \cite{T} p.\,371).
The best unconditional omega result for $M(x)$ is
\[
  M(x) = \Omega_{\pm}(x^{\frac{1}{2}}) \ .
\]
It should also be noted that if $\zeta(s)$ had a multiple zero of
of order $m \ge 2$ then
\[
    M(x) = \Omega_{\pm}(x^{\frac{1}{2}} (\log x)^{m-1}) \ .
\]
However, if RH is false then
\[
   M(x) = \Omega_{\pm}(x^{\theta - \delta})
\]
where
\[
  \theta = \sup_{\rho , \zeta(\rho) = 0} \mathrm{Re}(\rho)
\]
and $\delta$ is any positive constant (see Ingham \cite{In0} p.\
90).

To better understand the behaviour of $M(x)$, it is useful to
consider the closely related function
\begin{equation}
  \psi(x) - x = \sum_{n \le x} \Lambda(n) - x
\end{equation}
where $\Lambda(n)$ is Von-Mangoldt's function defined by
\begin{equation}
 \Lambda(n) = \left\{ \begin{array}{ll}
                  \log p & \mbox{if $n=p^{j} \ , \ j \ge 0$} \\
                  0 & \mbox{otherwise} \\
                  \end{array}
          \right.   \ .
\end{equation}
Here we review what is known concerning $\psi(x)-x$.
This may give us some better idea what type of upper and lower
bounds we should expect for $M(x)$.
Von Koch (see \cite{D} p.116) showed that RH is equivalent to
\begin{equation}
   \psi(x) - x =  O(x^{\frac{1}{2}} \log^{2} x ) \ .
\end{equation}
Moreover, Gallagher \cite{Ga} showed that RH implies
that
\[ \psi(x) - x =  O ( x^{\frac{1}{2}} (\log \log x)^{2} )
  \]
except on a set of finite logarithmic measure.
On the other hand, Littlewood demonstrated
\[ \psi(x) - x = \Omega_{\pm}
   \left( x^{\frac{1}{2}} \log \log \log x \right)
\]
under the assumption of RH (see \cite{In0} Chapter V) .  Moreover,
Montgomery \cite{Mo} has given an unpublished
probabilistic argument that suggests
\begin{equation}
  \underline{\overline{\lim}}
  \frac{\psi(x) - x}{x^{\frac{1}{2}} (\log \log \log x)^{2}}
  = \pm \frac{1}{2 \pi}
  \label{eq:mc}
\end{equation}
under the assumption of the Riemann hypothesis and the LI
conjecture.

Although the Mertens conjecture
is false, we can still obtain some averaged upper bounds
for $M(x)$.  We prove the following results:
\newtheorem{weakmertens}{Theorem}
\begin{weakmertens}
The Riemann Hypothesis and $J_{-1}(T) \ll T$ imply: \\
$(i)$
\[
   M(x) \ll x^{\frac{1}{2}} (\log x)^{\frac{3}{2}} \ ,
\]
$(ii)$
\[
   M(x) \ll x^{\frac{1}{2}} (\log \log x)^{\frac{3}{2}}
\]
except on a set of finite logarithmic measure, \\
(iii)
\[ \int_{2}^{X} \frac{M(x)^{2}}{x} \ dx \ll X  \ ,
\]
(iv) and the weak Mertens conjecture~(\ref{eq:wmc})
\[
   \int_{2}^{X} \left( \frac{M(x)}{x} \right)^{2} \ dx
   \ll \log X \ .
   \label{eq:ca}
\]
\end{weakmertens}
Theorem 1$(i)$ is due to Gonek, but had never published. The proof
of Theorem 1 $(ii)$ follows an argument due to Gallagher \cite{Ga}
and the proofs of Theorem 1 $(iii),(iv)$ follow an argument due to
Cram\'{e}r \cite{Cr}.  We note that by a more careful calculation
we can obtain an asymptotic evaluation in $(iv)$. However, since
$(iv)$ is easily deduced from Lemma 6, we include the argument.

Our study of  $M(x)$ requires some notions from
probability theory. Most importantly, we make use of distribution functions.
A distribution function $F(x)$ on $\mathbb{R}$
satisfies $F$ is non-decreasing,
$ F(-\infty) = 0$, $F(\infty) = 1$, $F$ is right-continuous, and
 $F$ has a limit on the left at each $x \in \mathbb{R}$.
Recall that if $P$ is a probability measure on $\mathbb{R}$, then
$F_{P}(x) := P((-\infty,x])$ is a distribution function. On the
other hand,  given a distribution function $F(x)$, there is a
theorem from probability theory which states there exists a
probability measure $P$ on $\mathbb{R}$ such that $F=F_{P}$.

In an attempt to better
understand $M(x)$, we give a conditional
proof of the existence of a limiting distribution function for
$\phi(y) = e^{-\frac{y}{2}}M(e^{y})$.
The idea to prove such a theorem originated with Heath-Brown's
comment:
\cite{HB}
``It appears to be an open question whether
\[ x^{-\frac{1}{2}}M(x) = x^{-\frac{1}{2}}\sum_{n \le x} \mu(n)
\]
has a distribution function.  To prove this one would want to
assume the Riemann Hypothesis and the simplicity of the zeros, and
perhaps also a growth condition on $M(x)$.'' Applying techniques
from Cram\'{e}r \cite{Cr} and Rubinstein-Sarnak \cite{RS} we
establish the following result.
\newtheorem{limdis}[weakmertens]{Theorem}
\begin{limdis}
Assume the Riemann Hypothesis and $J_{-1}(T) \ll T$.
Then $e^{-\frac{y}{2}}M(e^{y})$ has a limiting
distribution $\nu$ on $\mathbb{R}$, that is,
\begin{equation}
   \lim_{Y \rightarrow \infty} \frac{1}{Y} \int_{0}^{Y}
   f(e^{-\frac{y}{2}}M(e^{y})) \ dy =
   \int_{-\infty}^{\infty} f(x) \ d\nu(x)
   \label{eq:ng}
\end{equation}
for all bounded Lipschitz continuous functions $f$ on $\mathbb{R}$.
\end{limdis}
We note that the above theorem may be extended to all bounded
continuous functions $f(x)$ by standard approximation techniques.
However, we omit these arguments to keep the exposition simple.
Clearly Theorem 2 is useful in studying $M(x)$.
To see this,
suppose the above theorem remains valid for indicator functions.
Let $V$ be a fixed real number
and define $f = 1_{V}$ where
\[ 1_{V}(x) =
   \left\{  \begin{array}{ll}
          1 \ \mathrm{if} \ x \ge V \\
          0 \ \mathrm{if} \ x < V
          \end{array}  \right. \ .
\]
With the above choice of $f(x)$~(\ref{eq:ng}) translates to
\begin{equation}
  \lim_{Y \rightarrow \infty} \frac{1}
   {Y} \mathrm{meas} \{ y \in [0 , Y] \ | \
   M(e^{y}) \ge e^{\frac{y}{2}}V \ \}  =
   \nu([V,\infty)) .
   \label{eq:mcv}
\end{equation}
As noted in \cite{RS} p.\,174, the above identity would be true if
$\nu(x)$ is absolutely continuous. Under the additional assumption
of LI, one may show that $\nu$ is absolutely continuous. Moreover,
the LI conjecture implies that the Fourier transform of $\nu$ may
be computed explicitly.
\newtheorem{limdis2}{Corollary}
\begin{limdis2}
Assume the Riemann Hypothesis, $J_{-1}(T) \ll T$, and LI.
Then the Fourier
transform $\widehat{\nu}(\xi) = \int_{-\infty}^{\infty} e^{- i \xi t}
  \ d\nu(t)$
exists and equals
\begin{equation}
   \widehat{\nu}(\xi) = \prod_{\gamma > 0} \tilde{J}_{0} \left( \frac{2 \xi}
   {\left| (\frac{1}{2}+i\gamma)  \zeta^{'}(\cc+i\gamma) \right|}
   \right)
   \label{eq:ipro}
\end{equation}
where $\tilde{J}_{0}(z)$ is the Bessel function
$\tilde{J}_{0}(z) = \sum_{m=0}^{\infty} \frac{(-1)^{m}(\frac{1}{2}z)^{2m}}
                                {(m!)^{2}}$.
\end{limdis2}
Note that we have employed non-standard notation for the Bessel
function, so as not to confuse it with the moments $J_{-k}(T)$.
Under the same assumptions as Corollary 1, we observe that the set
\[ S = \{ x \ge 1 \ | \ |M(x)| \le \sqrt{x} \ \}
\]
has a logarithmic density.  Namely,
\[
  \delta(S) = \lim_{X \to \infty} \frac{1}{\log X}
  \int_{[2,X] \bigcap S} \frac{dt}{t}
\]
exists and $0 < \delta(S) < 1$. Since no counterexamples to the
Mertens conjecture have ever been found, we expect this
logarithmic density to be very close to $1$. In fact, preliminary
calculations indicate this.

In the same spirit as Theorems 1 and 2 we prove a strong
form of the weak Mertens conjecture is true.
This follows Cram\'{e}r's argument \cite{Cr}
subject to the same assumptions as the previous theorems.
\newtheorem{wm}[weakmertens]{Theorem}
\begin{wm}
Assume the Riemann hypothesis and $J_{-1}(T) \ll T$, then we have
\begin{equation}
   \int_{0}^{Y} \left( \frac{M(e^{y})}{e^{\frac{y}{2}}} \right)^{2}
   \, dy \sim \beta Y
   \label{eq:sqr}
\end{equation}
where
\begin{equation}
   \beta = \sum_{\gamma > 0} \frac{2}{|\rho \zeta^{'}(\rho)|^{2}}
   \ .
   \label{eq:zsum}
\end{equation}
Note that the assumption $J_{-1}(T) \ll T$ implies~(\ref{eq:zsum})
is convergent.
\end{wm}
A change of variable transforms~(\ref{eq:sqr}) to
\begin{equation}
   \int_{1}^{X} \left( \frac{M(x)}{x} \right)^{2} \,
   dx \sim
   \beta \log X  \ .
\end{equation}
Also, note that Theorem 3 corresponds to Theorem 2 with $f(x) =
x^{2}$. However, $f(x) = x^{2}$ is not a bounded function and does
not fall under the assumptions of Theorem 2.
We further note that the same techniques allow one to
establish
\begin{equation}
    \int_{0}^{Y} \frac{M(e^{y})}{e^{\frac{y}{2}}}
   \, dy = o(Y)
   \label{eq:fir}
\end{equation}
under the same conditions as Theorem 3.
Consequently,~(\ref{eq:sqr}) and~(\ref{eq:fir}) reveal that the
variance of the probability measure constructed in Theorem 2 is
$\beta$.

As one can see by  equation~(\ref{eq:mcv}) and Theorem 3,
the constructed limiting distribution of Theorem 2 reveals significant
information concerning $M(x)$.  The above formula~(\ref{eq:ipro})
for the Fourier
transform is crucial in studying the behaviour of $x^{-\frac{1}{2}}M(x)$.
Upon proving Theorem 2, we realized that the
constructed distribution could be used to study
large values of $M(x)$.
Using Montgomery's probabilistic methods we study
the tail of this distribution and give a conditional proof that
\[
   \exp(-\exp(\tilde{c}_{1} V^{\frac{4}{5}}))
   \ll \nu ( [V, \infty))
   \ll
   \exp(-\exp(\tilde{c}_{2} V^{\frac{4}{5}}))
\]
for positive effective constants $\tilde{c}_{1}$ and
$\tilde{c}_{2}$. We believe that $\tilde{c}_{1} = \tilde{c}_{2}$,
however it is not presently clear what this value should be.
Nevertheless, these bounds seem to suggest the following version
of an unpublished conjecture of Gonek's.

\noindent \textbf{Gonek's Conjecture} There exists a number $B >
0$ such that
\begin{equation}
    \underline{\overline{\lim}}_{x \to \infty} \frac{M(x)}{\sqrt{x}
   (\log \log \log
    x)^{\frac{5}{4}}} = \pm B \ .
    \label{eq:lb}
\end{equation}
After the completion of this work
the author learned that Gonek had arrived at this conjecture
at least ten years ago via Montgomery's techniques. He
had annuciated this conjecture at several conferences
in the early 1990's.
We note that the
exponent of the iterated triple logarithm is $\frac{5}{4}$
in~(\ref{eq:lb}) precisely because of the
Gonek-Hejhal conjecture~(\ref{eq:gh}).
Montgomery's conjecture~(\ref{eq:mc})
on $\frac{\psi(x)-x}{\sqrt{x}}$ shows that
the corresponding exponent on the iterated triple logarithm is $2$.
The difference between these cases is due directly to the different
discrete moments of
\[
  \sum_{\gamma \le T} \frac{1}{|\rho|} \asymp (\log T)^{2} \ \mathrm{and}
  \
  \sum_{\gamma \le T} \frac{1}{|\rho \zeta^{'}(\rho)|}
  \asymp  (\log T)^{\frac{5}{4}}
\]
where the second inequality is currently conjectural.

Finally, we remark that many of the results in this paper may be
extended to the summatory function of the Liouville function.
The Liouville function is defined as
$\lambda(n) = (-1)^{\Omega(n)}$ where $\Omega(n)$
denotes the total number of prime factors of $n$. P\'{o}lya was interested
in the summatory function
\begin{equation}
   L(x) = \sum_{n \le x} \lambda(n)
\end{equation}
since if the inequality $L(x) \le 0$ always persisted then the Riemann
hypothesis would follow.  Haselgrove \cite{Ha}
showed that this statement cannot
be true.  By the methods of this article, we can prove that
$e^{-\frac{y}{2}} L(e^{y})$ has a limiting distribution under the same
conditions as Theorems 1-3.  The reason we can extend the work to
this case is because
\begin{equation}
     \frac{\zeta(2s)}{\zeta(s)} = \sum_{n=1}^{\infty}
     \frac{\lambda(n)}{n^{s}}
\end{equation}
and thus the only difference is the term $\zeta(2s)$ in the
numerator. Nevertheless, this can be treated easily since we
understand the zeta function on the $\mathrm{Re}(s)=1$ line.

\emph{The majority of this article constitutes the last chapter of
my Ph.D. thesis. However, Theorem 3 was proven during a stay at
the Institute for Advanced Study during Spring 2002. I would like
to thank my Ph.D. supervisor, Professor David Boyd, for providing
me with academic and financial support during the writing of the
thesis and throughout my graduate studies. I also thank the I.A.S.
for its support and excellent working conditions. Finally, thanks
to Professor Steve Gonek for allowing me include Theorem 1$(i)$
and also for informing me of his earlier unpublished work.}

\section{Proof of Theorem 1}

Various estimates throughout this article require estimates for
averages of sums containing the expression $|\zeta^{'}(\rho)|^{-1}$.
This lemma establishes such estimates, subject to various special
cases of the Gonek-Hejhal conjecture~(\ref{eq:gh}).
\newtheorem{asym}{Lemma}
\begin{asym}
$(i)$ $J_{-\frac{1}{2}}(T) = \sum_{0 < \gamma < T}
    |\zeta^{'}(\rho)|^{-1} \ll
    T(\log T)^{v}$ implies
\[ \sum_{0 < \gamma < T} \frac{1}{|\rho \zeta^{'}(\rho)|} \ll
   (\log T)^{v+1} .
\]
$(ii)$ $J_{-1}(T) = \sum_{0 < \gamma < T}
     |\zeta^{'}(\rho)|^{-2} \ll T$ implies
\[  \sum_{T < \gamma < 2T} \frac{1}{|\rho \zeta^{'}(\rho)|^{2}}
    \ll \frac{1}{T} .
\]
$(iii)$ $J_{-\frac{1}{2}}(T) = \sum_{0 < \gamma < T}
      |\zeta^{'}(\rho)|^{-1} \ll T^{u}(\log T)^{v} $
    implies
\[
   \sum_{\gamma > T} \frac{(\log \gamma)^{a}}{\gamma^{b}
   |\zeta^{'}(\rho)|} \ll
   \frac{(\log T)^{a + v}}{T^{b-u}}
\]
subject to  $b > u \ge 1$.
\end{asym}
\pr For part$(i)$ note that
\begin{equation}
\begin{split}
   \sum_{0 < \gamma < T} \frac{1}{|\rho \zeta^{'}(\rho)|} \ll
   \sum_{0 < \gamma < T} \frac{1}{|\zeta^{'}(\rho)| \gamma}
   & = \left[ \frac{J_{-\frac{1}{2}}(t)}{t}
               \right]_{14}^{T}
   + \int_{14}^{T} \frac{J_{-\frac{1}{2}}(t)}{t^{2}} \ dt \\
   & = O \left( (\log T)^{v} +
   \int_{14}^{T} \frac{t(\log t)^{v}
   }{t^{2}} \ dt \right) \\
   & = O( (\log T)^{v+1} ) \ . \\
\end{split}
\end{equation}
Observe that we have made use of fact that all non-trivial zeros
$\rho = \beta + i \gamma$ satisfy $|\gamma| \ge 14$.
Part $(ii)$ is proven in an analogous fashion. For part
$(iii)$ let $\phi(t) = (\log t)^{a} t^{-b}$ and note that
its derivative is
$\phi^{'}(t) = (a(\log t)^{a-1} - b(\log t)^{a})/t^{b+1}$ .
Partial summation implies
\[ \sum_{\gamma > T} \frac{(\log \gamma)^{a}}{\gamma^{b}
   |\zeta^{'}(\rho)|}
   = \left[ \phi(t)J_{-\frac{1}{2}}(t) \right]_{T}^{\infty}
   - \int_{T}^{\infty} J_{-\frac{1}{2}}(t) \phi^{'}(t) \ dt .
\]
The first term is $\ll \phi(T) J_{-\frac{1}{2}}(T) = (\log
T)^{a+v}/T^{b-u}$.  Assuming the bound on $J_{-\frac{1}{2}}(T)$,
the second term is
\[ \ll \int_{T}^{\infty}
   \frac{ (t^{u}(\log t)^{v}) (\log t)^{a}}{t^{b+1}} \ dt
   = \int_{T}^{\infty} \frac{ (\log t)^{a+v}}{t^{b-u+1}} \ dt
   \ll \frac{(\log T)^{a+v}}{T^{b-u}}
\]
where the last integral is computed by an integration by parts.

We require Perron's formula in order to express $M(x)$ as the sum
of a complex integral and an error term.
\newtheorem{Perron}[asym]{Lemma}
\begin{Perron}
Let $f(s) = \sum_{n =1}^{\infty} a_{n} n^{-s}$ be absolutely
convergent for $\sigma = \mathrm{Re}(s) > 1$, $a_{n} \ll \Phi(n)$
where $\Phi(x)$ is positive and non-decreasing, and
\[ \sum_{n=1}^{\infty} \frac{|a_{n}|}{n^{\sigma}} =
   O \left( \frac{1}{(\sigma-1)^{\alpha}} \right) \ \mathrm{as} \
   \sigma \to 1^{+} \ .
\]
Then if $w = u+iv$ with $c > 0$, $u+c >1$, $T > 0$, we have for
all $x \ge 1$
\begin{equation}
\begin{split}
   \sum_{n \le x} \frac{a_{n}}{n^{w}} & =
   \frac{1}{2 \pi i} \int_{c-iT}^{c+iT} f(w+s) \frac{x^{s}}{s} ds \\
   & + O \left( \frac{x^{c}}{T(u+c-1)^{\alpha}} +
   \frac{\Phi(2x)x^{1-u} \log(2x)}{T} +
   \Phi(2x)x^{-u} \right) \\
\end{split}
\end{equation}
\end{Perron}
\pr This is a well-known theorem and is proven in \cite{Pr}
pp.\,376-379.

We need the following technical lemma in order to choose a good
contour for the complex integral obtained by Perron's formula.
\newtheorem{tn}[asym]{Lemma}
\begin{tn}
There exists a sequence of numbers
$\mathcal{T} = \{ T_{n} \}_{n=0}^{\infty}$ which satisfies
\[  n \le T_{n} \le n+1 \  and  \
   \frac{1}{\zeta(\sigma + i T)} = O(T^{\epsilon})
\]
for all $-1 \le \sigma \le 2$.
\end{tn}
\pr The above fact is proven in Titchmarsh \cite{T} pp.\,357-358
in the range $\frac{1}{2} \le \sigma \le 2$. It remains to prove
the bound in the range $-1 \le \sigma < \frac{1}{2}$. The
asymmetric form of the functional equation of the zeta function is
$\zeta(s) = \chi(s) \zeta(1-s)$ where $\chi(s) =
\pi^{s-\frac{1}{2}} \Gamma(\frac{1-s}{2}) / \Gamma(\frac{s}{2})$.
A calculation with Stirling's formula demonstrates that
$|\chi(\sigma+iT)| \asymp T^{\frac{1}{2} - \sigma}$ and therefore
we deduce that
\[
  |\zeta(s)|^{-1} = |\zeta(1-s) \chi(s)|^{-1}
  \ll T^{\epsilon + \sigma-\frac{1}{2}}
  \ll T^{\epsilon}
\]
for  $-1 \le \sigma < \frac{1}{2}$.

We now prove an explicit formula for $M(x)$.  With the exception
of a few minor changes, the proof follows Theorem 14.27 of
\cite{T} pp.\,372-374.
\newtheorem{Mertens}[asym]{Lemma}
\begin{Mertens}
Assume the Riemann hypothesis and that all zeros of
$\zeta(s)$ are simple. For $x \ge 2$ and $T \in \mathcal{T}$
\[ M(x) = \sum_{|\gamma| < T} \frac{x^{\rho}}{\rho \zeta^{'}(\rho)}
   + \tilde{E}(x,T)
\]
where
\[ \tilde{E}(x,T) \ll
   \frac{x \log x}{T} + \frac{x}{T^{1 - \epsilon}\log x}
   + 1  \ .
\]
\end{Mertens}
\pr  We apply Lemma 2 with $f(s) = \zeta(s)^{-1}$, $\Phi(x)=1$,
$\alpha=1$, and $w=0$ to obtain
\[ M(x) = \frac{1}{2 \pi i} \int_{c - iT}^{c + iT} \frac{x^{s}}{s \zeta(s)} ds
   + O \left( \frac{x^{c}}{T(c -1)}
   + \frac{x \log x}{T}
   + 1 \right) .
\]
Setting $c = 1 + (\log x)^{-1}$, this becomes
\[ M(x) = \frac{1}{2\pi i} \int_{c - iT}^{c + iT} \frac{x^{s}}{s \zeta(s)} ds
   + O \left( \frac{x \log x}{T} + 1 \right) .
\]
We introduce a large parameter $U$ and consider a positively
oriented rectangle $B_{T,U}$ with vertices at $c - iT,c + iT,-U +
iT,$ and $-U - i T$.  Thus the integral on the right equals
\[
    \frac{1}{2\pi i} \int_{B_{T,U}} \frac{x^{s}}{s \zeta(s)} ds
      - \frac{1}{2\pi i} \left( \int_{c + iT}^{-U + iT} +
                           \int_{-U + iT}^{-U - iT} +
                           \int_{-U - iT}^{c - iT} \right)
    \frac{x^{s}}{s \zeta(s)} ds  \ .
\]
It is shown in Titchmarsh \cite{T} p.\,373  that the middle
integral approaches 0 as $U \rightarrow \infty$. Inside the box
$B_{T,U}$, $\frac{x^{s}}{s \zeta(s)}$ has poles at the zeros of
the zeta function and $s=0$. By Cauchy's Residue Theorem, we have
\begin{equation}
\begin{split}
    M(x) & = \frac{1}{2\pi i} \sum_{|\gamma| < T}
    \frac{x^{\rho}}{\rho \zeta^{'}(\rho)}
    - 2 + \sum_{k \ge 1} \frac{x^{-2k}}{(-2k) \zeta^{'}(-2k)} \\
    &  - \frac{1}{2\pi i} \left( \int_{c + iT}^{-\infty + iT} +
                           \int_{-\infty - iT}^{c - iT} \right)
    \frac{x^{s}}{\zeta(s)s} ds
   + O \left( \frac{x \log x}{T} + 1 \right)  . \\
\end{split}
\end{equation}
The second and third terms are absorbed by the $O(1)$ term. We now
bound the integrals. Break up the first integral in two pieces as
\begin{equation}
  \int_{c+iT}^{-\infty + iT} \frac{x^{s}}{s \zeta(s)} \ ds =
   \left( \int_{c + iT}^{-1 + i T} +
   \int_{-1 + iT}^{-\infty + iT} \right) \frac{x^{s}}{s \zeta(s)} \ ds
   \ .
   \label{eq:2in}
\end{equation}
By Lemma 3, we have
\[ \left| \int_{-1 + i T}^{c + i T} \frac{x^{s}}{s \zeta(s)} \ ds
\right|
   \ll \int_{-1}^{c} \frac{x^{\sigma}T^{\epsilon}}
    { \sqrt{\sigma^{2} + T^{2}} } \ d \sigma \le T^{\epsilon -1}
   \int_{-1}^{c} e^{\sigma \log x} \ d \sigma \ll
    \frac{x}{T^{1-\epsilon}\log x} .
\]
For the second piece we apply the functional equation
\[
   \int_{-\infty + iT}^{-1 + iT} \frac{x^{s}}{s \zeta(s)} \ ds
    = \int_{2 - iT}^{\infty - i T}
    \frac{x^{1-s}2^{s-1} \pi^{s}}{(1-s)\cos(\frac{\pi s}{2}) \Gamma(s) \zeta(s)} \
     ds .
\]
For $\sigma \ge 2$ we have the Stirling formula estimate
$\frac{1}{|\Gamma(\sigma - iT)|} \ll e^{\sigma -
(\sigma-\frac{1}{2})\log \sigma + \frac{1}{2} \pi T}$ and the
elementary estimate $\frac{1}{|\cos(\frac{\pi (\sigma - iT)}{2})|}
\ll e^{-\frac{\pi}{2}T}$ and hence the integral is
\[ O \left( \int_{2}^{\infty} \frac{x}{T}
   \left( \frac{2 \pi}{x} \right)^{\sigma}
   e^{\sigma - (\sigma-\frac{1}{2})\log \sigma} \ d \sigma \right)
   = O \left( \frac{x}{T} \right) \ .
\]
The same argument applies to the second integral in~(\ref{eq:2in}) and we
have shown that
\[ M(x) = \sum_{|\gamma| < T} \frac{x^{\rho}}{\rho \zeta^{'}(\rho)}
   + O\left( \frac{x \log x}{T}
   + \frac{x}{T^{1 - \epsilon}\log x}
   + 1 \right) \ .
\]

We now remove the assumption that $T \in \mathcal{T}$ from the last
lemma by applying the Gonek-Hejhal conjecture~(\ref{eq:gh}) for $k=-1$.
\newtheorem{Mertens2}[asym]{Lemma}
\begin{Mertens2}
Assume the Riemann hypothesis and
\[ J_{-1}(T) \ll T \ .
\]
For $x \ge 2$, $T \ge 2$
\[ M(x) = \sum_{|\gamma| \le T} \frac{x^{\rho}}{\rho \zeta^{'}(\rho)}
   + E(x,T)
\]
where
\begin{equation}
   E(x,T) \ll
   \frac{x \log x}{T} + \frac{x}{T^{1 - \epsilon}\log x}
   + \left(\frac{x \log T}{T} \right)^{\frac{1}{2}} + 1 \ .
   \label{eq:erm}
\end{equation}
\end{Mertens2}
\pr Let $T \ge 2$ satisfy $n \le T \le n + 1$.  Now suppose
without loss of generality that $n \le T_{n} \le T \le
n + 1$. Then we have
\[ M(x) = \sum_{|\gamma| \le T} \frac{x^{\rho}}{\rho \zeta^{'}(\rho)}
  - \sum_{T_{n} \le |\gamma| \le T} \frac{x^{\rho}}{\rho \zeta^{'}(\rho)}
  + \widetilde{E}(x,T_{n}) \ .
\]
By Cauchy-Schwarz the second sum is
\[
   \left| \sum_{T_{n} \le \gamma \le T}
   \frac{x^{\rho}}{\rho \zeta^{'}(\rho)} \right|
   \le x^{\frac{1}{2}} \left(
   \sum_{T_{n}  \le \gamma \le T}
   \frac{1}{|\rho \zeta^{'}(\rho)|^{2}} \right)^{\frac{1}{2}}
   \left( \sum_{T_{n} \le \gamma \le T}  1 \right)^{\frac{1}{2}}
   \ .
\]
By Lemma 1$(ii)$, $J_{-1}(T) \ll T$ implies $\sum_{T \le \gamma
\le
  2T} \frac{1}{|\rho \zeta^{'}(\rho)|^{2}} \ll \frac{1}{T}$ and we
deduce that
\[ \left| \sum_{T_{\nu} \le \gamma \le T}
   \frac{x^{\rho}}{\rho \zeta^{'}(\rho)} \right| \ll
    \left( \frac{x \log T}{T} \right)^{\frac{1}{2}}
\]
which completes the proof.

Lemma 6 is the crucial step in proving the existence of the
limiting distribution in the next section.  The key point is that
the integral in this lemma should be small in order to justify the
weak convergence of a sequence of distribution functions in
Theorem 2. This is also used in the proof of Theorem 1 parts
$(ii)$-$(iv)$.
\newtheorem{ub}[asym]{Lemma}
\begin{ub}
Assume the Riemann hypothesis and $J_{-1}(T) \ll T$. Then
\begin{equation}
   \int_{Z}^{eZ} \left| \sum_{T \le |\gamma| \le X}
   \frac{x^{i \gamma}}{\rho \zeta^{'}(\rho)} \right|^{2}
   \frac{dx}{x} \ll
   \frac{(\log T)}{T^{\frac{1}{4}}}
   \label{eq:kl}
\end{equation}
for $Z > 0$ and $T < X$.
\end{ub}
\pr Making the substitution $x = e^{y}$ in the left hand side
of~(\ref{eq:kl})
we obtain
\begin{equation}
\begin{split}
   &  \int_{\log Z}^{\log Z + 1} \left| \sum_{T \le |\gamma| \le X}
   \frac{e^{i \gamma y}}{\rho \zeta^{'}(\rho)} \right|^{2} \ dy
   \le 4 \int_{\log Z}^{\log Z + 1} \left| \sum_{T \le \gamma \le X}
     \frac{e^{i \gamma y}}{\rho \zeta^{'}(\rho)} \right|^{2} \ dy \\
   & = 4 \sum_{T \le \gamma \le X} \sum_{T \le \lambda \le X}
     \frac{1}{\rho \zeta^{'}(\rho) \overline{\rho^{'}}
     \overline{\zeta^{'}(\rho^{'})} }
     \int_{\log Z}^{\log Z + 1} e^{i(\gamma - \lambda) y} \ dy  \\
  & \ll \sum_{T \le \gamma \le X} \sum_{T \le \lambda \le X}
     \frac{1}{|\rho \zeta^{'}(\rho)|| \overline{\rho^{'}}
     \overline{\zeta^{'}(\rho^{'})}| }
     \mathrm{min}  \left( 1,\frac{1}{|\gamma - \lambda|} \right) \\
\end{split}
\end{equation}
Note that $\rho$ and $\rho^{'}$ denote zeros of the form
$\rho = \frac{1}{2} + i \gamma$ and $\rho^{'} = \frac{1}{2} + i \lambda$.
We break this last sum in
two sums  $\Sigma_{1}$ and $\Sigma_{2}$ where
$\Sigma_{1}$ consists of those terms for which
$|\gamma - \lambda| \le 1$ and  $\Sigma_{2}$ consists of the
complementary
set.  The first sum is bounded as follows
\[ \Sigma_{1} \ll \sum_{T \le \gamma \le X}
   \frac{1}{|\rho \zeta^{'}(\rho)|}
   \sum_{\gamma - 1 \le \lambda \le \gamma + 1}
   \frac{1}{|\rho^{'} \zeta^{'}(\rho^{'})|} .
\]
It is well known that $N(t+1) - N(t-1) \ll \log t$, hence the
inner sum is
\[
   \le
   \left( \sum_{\gamma - 1 \le \lambda \le \gamma + 1}
   \frac{1}{|\rho^{'} \zeta^{'}(\rho^{'})|^{2}} \right)^{\frac{1}{2}}
   ( N(\gamma + 1) - N(\gamma - 1))^{\frac{1}{2}}
   \ll
   \left( \frac{\log \gamma}{\gamma} \right)^{\frac{1}{2}}
\]
by an application of Lemma 1$(ii)$. By Lemma 1$(iii)$ we deduce
that
\begin{equation}
   \Sigma_{1}  \ll   \sum_{T \le \gamma}
   \frac{(\log \gamma )^{\frac{1}{2}}}
   {\gamma^{\frac{3}{2}} |\zeta^{'}(\rho)|}
   \ll \frac{\log T}{T^{\frac{1}{2}}} \ .
\end{equation}
 Write the second sum as
\begin{equation}
\begin{split}
   \Sigma_{2}
   & =
   \sum_{T \le \gamma \le X} \frac{1}{|\rho\zeta^{'}(\rho)|}
   \sum_{T \le \lambda \le X, |\gamma - \lambda| \ge 1}
   \frac{1}
   {|\rho^{'} \zeta^{'}(\rho^{'})| | \gamma - \lambda |} .\\
   \label{eq:sig2}
\end{split}
\end{equation}
The inner sum is analyzed by splitting the sum in to ranges. The
crucial range is when $|\gamma-\lambda| \approx 1$. This argument
was originally employed by Cram\'{e}r \cite{Cr}. We eliminate the
condition $\gamma \le X$ and denote the inner sum
of~(\ref{eq:sig2}) as $S(\gamma)$ where $\gamma \ge T$. Consider
the set of numbers, $\gamma^{\cc}, \gamma - \gamma^{\cc}$, and
$\gamma - 1$. One of the following cases must occur $T \le
\gamma^{\cc}$, $\gamma^{\cc} < T \le \gamma - \gamma^{\cc}$, $
\gamma - \gamma^{\cc} < T \le \gamma - 1$, or $\gamma-1 < T \le
\gamma$.  These conditions translate in to the four cases: $T^{2}
\le \gamma$, $T + \sqrt{T+\frac{1}{4}} + \frac{1}{2} \le \gamma <
T^{2}$, $T+1 \le \gamma < T + \sqrt{T+\frac{1}{4}} + \frac{1}{2}$,
and $T \le \gamma < T+1$. Suppose the first case is true, i.e. $T
\le \gamma^{\cc}$. Then we may write the inner sum $S(\gamma)$ as
six separate sums
\begin{equation}
\begin{split}
   S(\gamma) = & \left(
   \sum_{T \le \lambda < \gamma^{\cc}} +
   \sum_{\gamma^{\cc} \le \lambda < \gamma - \gamma^{\cc}} +
     \sum_{\gamma - \gamma^{\cc}  \le \lambda \le \gamma - 1} \right.  \\
  &  \left. + \sum_{\gamma + 1 \le \lambda < \gamma +\gamma^{\cc}} +
   \sum_{\gamma +\gamma^{\cc} \le \lambda < 2 \gamma} +
   \sum_{2 \gamma \le \lambda} \right)
    \frac{1}
   {|\rho^{'} \zeta^{'}(\rho^{'})| | \gamma - \lambda |} \ . \\
\end{split}
\end{equation}
Denote these sums by $\sigma_{1}, \ldots , \sigma_{6}$. In the
following estimates we apply Lemma 1$(ii)$ several times. We find
that
\begin{equation}
\begin{split}
   \sigma_{1} \le \frac{1}{\gamma - \gamma^{\cc}}
   \sum_{T \le \lambda < \gamma^{\cc}}
    \frac{1}{|\rho^{'} \zeta^{'}(\rho^{'})|}
  & \ll \frac{1}{\gamma} \left(\sum_{T \le \lambda < \gamma^{\cc}}
   \frac{1}{|\rho^{'} \zeta^{'}(\rho)|^{2}} \right)^{\frac{1}{2}}
   \left( \sum_{T \le \lambda < \gamma^{\cc}} 1
   \right)^{\frac{1}{2}} \\
  &  \ll \frac{1}{\gamma T^{\frac{1}{2}}}
  (\gamma^{\cc} \log \gamma)^{\frac{1}{2}}
   = \frac{(\log \gamma)^{\frac{1}{2}}}{T^{\frac{1}{2}} \gamma^{\frac{3}{4}}} , \\
\end{split}
\end{equation}
\begin{equation}
\begin{split}
   \sigma_{2} \le \frac{1}{\gamma^{\cc}}
   \sum_{\gamma^{\cc} \le \lambda < \gamma - \gamma^{\cc}}
   \frac{1}{|\rho^{'} \zeta^{'}(\rho^{'})|}
  & \le  \frac{1}{\gamma^{\cc}}
   \left(\sum_{\gamma^{\cc} \le \lambda
   < \gamma - \gamma^{\cc}}\frac{1}{|\rho^{'}
   \zeta^{'}(\rho^{'})|^{2}} \right)^{\frac{1}{2}}
   \left( \sum_{\gamma^{\cc} \le \lambda <
   \gamma - \gamma^{\cc}} 1 \right)^{\frac{1}{2}} \\
   & \ll \frac{1}{\gamma^{\cc}} \left(
   \frac{1}{\gamma^{\cc}} \right)^{\frac{1}{2}}
   (\gamma \log \gamma)^{\frac{1}{2}} =
    \frac{(\log \gamma)^{\frac{1}{2}}}
   {\gamma^{\frac{1}{4}}} ,\\
\end{split}
\end{equation}
and
\begin{equation}
\begin{split}
   \sigma_{3} \le  \left(
   \sum_{\gamma - \gamma^{\cc} \le \lambda \le
   \gamma - 1}
    \frac{1}{|\rho^{'} \zeta^{'}(\rho^{'})|^{2}} \right)^{\frac{1}{2}}
   \left( \sum_{\gamma - \gamma^{\cc}
   \le \lambda \le \gamma - 1} 1 \right)^{\frac{1}{2}}
   \ll  \frac{1}{\gamma^{\frac{1}{2}}} (\gamma^{\cc}
   \log \gamma)^{\frac{1}{2}} =
   \frac{(\log \gamma)^\frac{1}{2}}{\gamma^{\frac{1}{4}}} .
\end{split}
\end{equation}
The fourth sum, $\sigma_{4}$, gives the same error as the third
sum. Similarly,
\[
   \sigma_{5}  \ll \frac{1}{\gamma^{\cc}}
   \left( \sum_{\gamma + \gamma^{\cc} \le \lambda}
   \frac{1}{|\rho^{'} \zeta^{'}(\rho)|^{2}} \right)^{\frac{1}{2}}
   \left( \sum_{\gamma + \gamma^{\cc} \le \lambda \le 2 \gamma}
   1 \right)^{\frac{1}{2}}
   \ll \frac{1}{\gamma^{\cc}}
   \left( \frac{\gamma \log \gamma}{\gamma} \right)^{\frac{1}{2}}
    = \frac{(\log \gamma)^{\frac{1}{2}}}{\gamma^{\cc}}
\]
and
\begin{equation}
\begin{split}
   \sigma_{6} & \le \sum_{k = 1}^{\infty}
   \sum_{2^{k}\gamma \le \lambda \le 2^{k+1}\gamma}
   \frac{1}{|\rho^{'} \zeta^{'}(\rho^{'})| | \gamma - \lambda |} \\
   & \le \sum_{k = 1}^{\infty} \frac{1}{(2^{k}-1)\gamma}
   \left( \sum_{2^{k}\gamma \le \lambda \le 2^{k+1}\gamma}
   \frac{1}{|\rho^{'} \zeta^{'}(\rho^{'})|^{2}} \right)^{\frac{1}{2}}
   \left( \sum_{2^{k}\gamma \le \lambda \le 2^{k+1}\gamma}
   1 \right)^{\frac{1}{2}} \\
   & \le \sum_{k = 1}^{\infty} \frac{1}{(2^{k}-1)\gamma}
   \left( \frac{2^{k+1}\gamma \log(2^{k+1}\gamma)}{2^{k}\gamma}
   \right)^{\frac{1}{2}} \ll
   \frac{(\log \gamma)^{\frac{1}{2}}}{\gamma} \ .  \\
\end{split}
\end{equation}
Putting together these bounds leads to
\[ S(\gamma) \ll \frac{(\log \gamma)^{\frac{1}{2}}}{\gamma^{\frac{1}{4}}}
\]
as long as $T^{2} \le \gamma$. In fact, the same argument applies
in the other three cases.  The only difference is that there would
be fewer sums and we still establish $S(\gamma) \ll (\log
\gamma)^{\frac{1}{2}} \gamma^{-\frac{1}{4}}$ for all $\gamma \ge
T$. The assumption $J_{-1}(T) \ll T$ implies by Cauchy-Schwarz
that
\begin{equation}
   J_{-1/2}(T) \ll J_{-1}(T)^{\frac{1}{2}} N(T)^{\frac{1}{2}}
   \ll T^{\frac{1}{2}} (T (\log T))^{\frac{1}{2}}
   = T (\log T)^{\frac{1}{2}} \ .
   \label{eq:jb}
\end{equation}
Applying Lemma 1$(iii)$ yields the bound
\[ \Sigma_{2} \ll \sum_{\gamma > T} \frac{(\log \gamma)^{\frac{1}{2}}}
   {\gamma^{\frac{5}{4}} |\zeta^{'}(\rho)|}
   \ll
   \frac{\log T}{T^{\frac{1}{4}}}
\]
and the lemma is proved.

Combining the previous lemmas we may now prove Theorem 1. \\
\noindent {\it Proof of Theorem 1.} $(i)$ By Lemma 5,
\[
  M(x) \ll x^{\frac{1}{2}} \sum_{0 < \gamma < T} \frac{1}{|\rho
  \zeta^{'}(\rho)|} + E(x,T)
\]
where $E(x,T)$ is defined by~(\ref{eq:erm}). By the
bound~(\ref{eq:jb}) Lemma 1$(i)$ yields
\[
   M(x) \ll x^{\frac{1}{2}} (\log T)^{\frac{3}{2}}
        + \frac{x \log x}{T} + \frac{x}{T^{1-\epsilon} \log x}
        + \left( \frac{x \log T}{T} \right)^{\frac{1}{2}} \ .
\]
By the choice $T^{1-\epsilon} = \sqrt{x}$, we deduce $M(x) \ll
\sqrt{x} (\log x)^{\frac{3}{2}}$. \\
\noindent $(ii)$ The starting point is to consider the explicit
formula. By Lemma 5, we have
\begin{equation}
   M(x) = \sum_{|\gamma| \le X}
   \frac{x^{\rho}}{\rho \zeta^{'}(\rho)} + O \left( X^{\epsilon}
   \right)
   \label{eq:ef}
\end{equation}
valid for $X \le x \ll X$.  By Lemma 6, we have
for $T^{4} < X$
\[
   \int_{X}^{eX} \left| \sum_{T^{4} \le |\gamma| \le X}
   \frac{x^{\rho}}{\rho \zeta^{'}(\rho)} \right|^{2}
   \frac{dx}{x^{2}} \ll
   \frac{(\log T)}{T} .
\]
By considering the set
\[
  S = \left\{  x \ge 2 \ | \
   \left| \sum_{T^{4} \le |\gamma| \le X}
   \frac{x^{\rho}}{\rho \zeta^{'}(\rho)} \right|  \ge
   x^{\frac{1}{2}}(\log \log x)^{\frac{5}{4}} \right\}
\]
it follows that
\[ (\log \log X)^{\frac{5}{2}} \int_{S \cap [X,eX]} \frac{dx}{x}
   \le \int_{X}^{eX} \left| \sum_{T^{4} \le |\gamma| \le X}
   \frac{x^{\rho}}{\rho \zeta^{'}(\rho)} \right|^{2}
   \frac{dx}{x^{2}} \ll
   \frac{(\log T)}{T}
\]
and thus
\[ \int_{S \cap [X,eX]} \frac{dx}{x} \ll
   \frac{ (\log T)}{T(\log \log X)^{\frac{5}{2}}}
   = \frac{1}{T (\log T)^{\frac{3}{2}}}
\]
for $T = \log X$.  Choosing $X = e^{k}$ with $k = 2,3, \ldots$ we
deduce
\[
   \int_{S \cap [e^{2},\infty]} \frac{dx}{x} \ll
   \sum_{k=2}^{\infty} \frac{1}{k (\log k)^{\frac{3}{2}}}
   < \infty
\]
and thus $S$ has finite logarithmic measure. By the
bound~(\ref{eq:jb}) Lemma 1$(i)$ implies
\[  \left| \sum_{0 \le |\gamma| \le T^{4}}
   \frac{x^{\rho}}{\rho \zeta^{'}(\rho)} \right| \ll
   X^{\frac{1}{2}}
   \sum_{0 \le |\gamma| \le T^{4}}
   \frac{1}{|\rho \zeta^{'}(\rho)|}
   \ll X^{\frac{1}{2}} (\log T)^{\frac{3}{2}}
   \ll X^{\frac{1}{2}} (\log \log X)^{\frac{3}{2}}
\]
for $X \le x \le eX$.  Hence,
\[ M(x) = \sum_{T^{4} \le |\gamma| \le X}
   \frac{x^{\rho}}{\rho \zeta^{'}(\rho)}
   + O \left( X^{\frac{1}{2}} (\log \log X)^{\frac{3}{2}}
   \right)
\]
for $X \le x \le eX$ and $T = \log X$.  Define the set
\[
  S_{\alpha} = \{  x \ge 2 \ | \
   |M(x)| \ge \alpha
    x^{\frac{1}{2}}(\log \log x)^{\frac{3}{2 }} \}.
\]
Suppose $x \in S_{\alpha} \bigcap
[X,eX]$.  Then we have
\begin{equation}
\begin{split}
  &  \left| \sum_{T^{4} \le |\gamma| \le X}
   \frac{x^{\rho}}{\rho \zeta^{'}(\rho)} \right|
    \ge  |M(x)|
   - O \left( X^{\frac{1}{2}} (\log \log X)^{\frac{3}{2}} \right)  \\
   &  \ge \alpha x^{\frac{1}{2}} (\log \log x)^{\frac{3}{2}}
   - O \left( X^{\frac{1}{2}} (\log \log X)^{\frac{3}{2}} \right)
   \ge  x^{\frac{1}{2}} (\log \log x)^{\frac{5}{4}} \\
   \label{eq:singe}
\end{split}
\end{equation}
for $x \in [X,eX]$ as long as $X$ is sufficiently large and
$\alpha$ is chosen larger than the constant that occurs in the
error term of~(\ref{eq:singe}). Thus $S_{\alpha} \bigcap [X,eX]
\subset S \bigcap [X,eX]$ for $X$ sufficiently large and it
follows that $S_{\alpha}$ has finite logarithmic measure.  Observe
that if we also assumed the conjecture $J_{-\frac{1}{2}}(t) \ll
t(\log t)^{\frac{1}{4}}$ then the same arguments in $(i)$ and
$(ii)$ would have shown that $M(x) \ll x^{\frac{1}{2}} (\log
x)^{\frac{5}{4}}$ and $M(x) \ll x^{\frac{1}{2}} (\log \log
x)^{\frac{5}{4}}$ except on a set of finite logarithmic measure.
\\
\noindent $(iii)$ Squaring equation~(\ref{eq:ef}), dividing by
$x^{2}$, and integrating yields
\begin{equation}
   \int_{X}^{eX} \left( \frac{M(x)}{x} \right)^{2} \ dx
   \ll \int_{X}^{eX} \left| \sum_{|\gamma| \le X}
   \frac{x^{\rho}}{\rho \zeta^{'}(\rho)} \right|^{2} \frac{dx}{x^{2}}
   + O(X^{-1+2\epsilon})
   \ll 1 \\
   \label{eq:cafe}
\end{equation}
by taking $Z=X$ and $T=14$ in Lemma 6.  It immediately follows that
\[ \int_{X}^{eX} \frac{M(x)^{2}}{x} \ dx
   \ll X .
\]
Substituting the values $\frac{X}{e}, \frac{X}{e^{2}}, \ldots$
and adding yields
\[
  \int_{2}^{X} \frac{M(x)^{2}}{x} dx \ll X \ .
\]
$(iv)$ Similarly, we obtain from~(\ref{eq:cafe})
\[
  \int_{2}^{X} \left( \frac{M(x)}{x} \right)^{2} dx
  \ll
  \sum_{k=1}^{[\log(\frac{X}{2})]+1} \int_{X/e^{k}}^{X/e^{k-1}}
  \left( \frac{M(x)}{x} \right)^{2}  dx \ll \log X \ .
\]
\section{Proofs of Theorems 2 and 3}

In this section
we prove the existence of a limiting distribution
for the function $\phi(y) = e^{-\frac{y}{2}}M(e^{y})$.
If we assume the Riemann hypothesis and write non-trivial zeros as $\rho = \frac{1}{2}
+ i \gamma$, then we obtain
\[ x^{-\frac{1}{2}}M(x) =
  \sum_{|\gamma| \le T} \frac{x^{i \gamma}}{\rho \zeta(\rho)}
  +  \sum_{T < |\gamma| \le e^{Y}} \frac{x^{i \gamma}}{\rho \zeta(\rho)}
  + x^{-\frac{1}{2}}E(x,e^{Y})
\]
where $T < e^{Y}$ and $E(x,e^{Y})$ is defined in~(\ref{eq:erm}).
Making the variable change $x = e^{y}$, we have
\begin{equation}
   \phi(y) = e^{-\frac{y}{2}}M(e^{y}) =
    \phi^{(T)}(y) + \epsilon^{(T)}(y)
\end{equation}
where
\begin{equation}
  \phi^{(T)}(y) = \sum_{|\gamma| \le T}
   \frac{e^{i\gamma y}}{\rho \zeta^{'}(\rho)} \ \  and \
   \label{eq:phity}
\end{equation}
\begin{equation}
  \epsilon^{(T)}(y) = \sum_{T \le |\gamma| \le e^{Y}}
   \frac{e^{i\gamma y}}{\rho \zeta^{'}(\rho)} +
   e^{-\frac{y}{2}}E(e^{y},e^{Y}) .
   \label{eq:ety}
\end{equation}

In order to construct a sequence of distribution functions that
converge to the distribution of Theorem 2, we require the
following uniform distribution result.
\newtheorem{kw}[asym]{Lemma}
\begin{kw}
Let $t_{1}, \ldots , t_{N}$ be $N$ arbitrary real numbers.
Consider the curve $\psi(y) = y(t_{1}, \ldots , t_{N}) \in
\mathbb{R}^{N}$ for $y \in \mathbb{R}$.
Let $f :\mathbb{R}^{N} \to \mathbb{R}$ be continuous and
have period one in each of its
variables.  There exists an integer $1 \le J \le N$ and
$A$, a $J$-dimensional parallelotope, such that
\begin{equation}
\begin{split}
   \lim_{Y \to \infty} \frac{1}{Y}
   \int_{0}^{Y} f(\psi(y)) \, dy
   & = \int_{A} f(a) \, d \mu  \\
\end{split}
\end{equation}
where $\mu$ is normalized Haar measure on $A$.  More
precisely, $A$ is the topological closure of $\psi(y)$ in
$\mathbb{T}^{N}$.
\end{kw}
{\it Proof}.  This lemma is a well-known and it is a variant of
the traditional Kronecker-Weyl theorem (see Hlawka \cite{Hl},
pp.\,1-14 for the proof) . We now describe the principal idea in
how the lemma is deduced from this.  Let $J$ be the maximum number
of linearly independent elements over $\Bbb Q$ among $t_{1},
\ldots , t_{N}$. The basic idea is to show that the topological
closure of the set $\{ \,
   \left( \{ y \frac{\gamma_{1}}{2 \pi} \} ,
   \ldots , \{ y \frac{\gamma_{N}}{2 \pi} \}  \right)
 \ | \ y \in \mathbb{R} \ \}$ cuts out a sub-torus of $\Bbb T^{N}$ of dimension
$J$ (Note that $\{x \}$ is the fractional part of $x \in \Bbb R$).
By a variable change, one then deduces the lemma from the
Kronecker-Weyl theorem.

By an application of Lemma 7, we construct for each large $T$ a
distribution function $\nu_{T}$.
\newtheorem{mer3}[asym]{Lemma}
\begin{mer3}
Assume the Riemann hypothesis, then
for each $T \ge \gamma_{1}$ (the imaginary ordinate of
the first non-trivial zero of $\zeta(s)$) there is a probability measure
$\nu_{T}$ on $\mathbb{R}$ such that
\[ \nu_{T}(f) := \int_{-\infty}^{\infty} f(x) \ d\nu_{T}(x)
   = \lim_{Y \rightarrow \infty} \frac{1}{Y} \int_{0}^{Y}
   f(\phi^{(T)}(y)) \ dy
\]
for all bounded continuous functions $f$ on $\mathbb{R}$ where
$\phi^{(T)}(y)$ is defined by~(\ref{eq:phity}).
\end{mer3}
\pr This is identical to Lemma 2.3 of \cite{RS} p.\,180. Let $N =
N(T)$ denote the number of zeros of $\zeta(s)$ to height $T$.
Label the imaginary ordinates of the zeros as $\{ \gamma_{1},
\ldots , \gamma_{N} \}$. By pairing conjugate zeros $\rho =
\frac{1}{2} + i \gamma$ and $\overline{\rho} = \frac{1}{2} - i
\gamma$ we have
\[ \phi^{(T)}(y) = \sum_{|\gamma| \le T}
   \frac{e^{i\gamma y}}{\rho \zeta^{'}(\rho)} =
   2 \mathrm{Re} \left(
   \sum_{l = 1}^{N} b_{l}
   e^{i y \gamma_{l}} \right)
\]
where  $b_{l} = \frac{1}{(\frac{1}{2} + i \gamma_{l})
\zeta^{'}(\frac{1}{2} + i \gamma_{l} )}$.
Define functions $X_{T}$ and $g$ on the
$N$-torus $\mathbb{T}^{N}$ by
\[ X_{T}(\theta_{1}, \ldots , \theta_{N}) =
   2 \mathrm{Re} \left(
   \sum_{l = 1}^{N} b_{l}
   e^{2 \pi i \theta_{l}} \right) \ \mathrm{and} \
   g(\theta_{1}, \ldots , \theta_{N}) = f(X_{T}(\theta_{1}, \ldots ,
   \theta_{N}))
   \ .
\]
We now apply Lemma 7 to the $N$ numbers $\{ \frac{\gamma_{1}}{2
\pi}, \ldots , \frac{\gamma_{N}}{2 \pi} \}$ and to the continuous
function $g$. According to Lemma 7 there exists a torus $A \subset
\mathbb{T}^{N}$ such that
\[
   \lim_{Y \to \infty} \frac{1}{Y} \int_{0}^{Y}
   g \left(y \left( \frac{\gamma_{1}}{2 \pi},
    \ldots , \frac{\gamma_{N}}{2 \pi} \right) \right) \, dy
   =  \int_{A} g(a) \, d \mu \ .
\]
The measure $d\mu$ is normalized Haar
measure on $A$.  Note that
$ \left.  X_{T} \right|_{A} : A \to \mathbb{R}$
is a random variable and we define a probability measure $\nu_{T}$ on
$\mathbb{R}$ by
\begin{equation}
   \nu_{T}(B) =  \mu \left( \left. X \right|_{A}^{-1} (B) \right)
   \label{eq:pullback}
\end{equation}
where $B$ is any Borel set.
By the change of variable formula, we deduce
\[
 \lim_{Y \rightarrow \infty} \frac{1}{Y} \int_{0}^{Y}
   f(\phi^{(T)}(y)) \, dy = \int_{\mathbb{R}} f(x) \, d \nu_{T}(x)
\]
and the proof is complete.

Before proceeding, we require
some results from probability theory.  We say that  a real valued
function $G(x)$ is
a generalized distribution function on $\mathbb{R}$ if it is
non-decreasing and right-continuous.  Lemma 9$(i)$
will enable us to construct a limiting distribution function
from the set $\{ \nu_{T} \}_{T \gg 1}$ constructed in the previous lemma.
\newtheorem{Prob}[asym]{Lemma}
\begin{Prob}
$(i)$ Let $\, F_{n} \,$ be a sequence of distribution functions.
There exists a subsequence $\{ F_{n_{k}} \}$  and a generalized
distribution function $F$ such that
\[ \lim_{k \to \infty} F_{n_{k}}(x) = F(x)
\]
at continuity points $x$ of $F$.  \\
$(ii)$ Let $\{ F_{n} \}$ be distribution functions and $F$ a
generalized distribution function on $\mathbb{R}$ such that
$F_{n}$ converges to $F$ weakly.  This is equivalent to
\[ \int_{\mathbb{R}} f(x) \, dF_{n}(x) \to \int_{\mathbb{R}}
   f(x) \, dF(x)
\]
for all continuous, bounded, real $f(x)$.  \\
$(iii)$ Let $F_{n}$, $F$ be distribution functions with Fourier transforms,
$\hat{F}_{n}$, $\hat{F}$.  A necessary and sufficient condition for $F_{n}$
to converge weakly to $F$ is $\hat{F}_{n}(t) \to \hat{F}(t)$ for each $t$.
\end{Prob}
\pr Part $(i)$ is Helly's selection theorem and part $(iii)$ is
Levy's theorem. See \cite{Bi} pp.\,344-346 for proofs of $(i)$ and
$(ii)$ and pp.\,359-360 for $(iii)$.

The next lemma shows that the error term
 $\epsilon^{(T)}(y)$ in~(\ref{eq:ety})
has small mean square.  This will be crucial in deducing that a limiting
distribution exists for $e^{-\frac{y}{2}}M(e^{y})$.
\newtheorem{mer2}[asym]{Lemma}
\begin{mer2}
Assume the Riemann hypothesis and $J_{-1}(T) \ll T$.
For $T \ge 2$ and $Y \ge \log 2$,
\[ \int_{\log 2}^{Y} |\epsilon^{(T)}(y)|^{2} dy \ll Y
   \frac{(\log T)}{T^{\frac{1}{4}}} + 1 .
\]
\end{mer2}
\pr First we will consider the contribution from $E(x,T)$ as
defined in~(\ref{eq:erm}). Note that
\[
 \int_{\log 2}^{Y} |e^{-\frac{y}{2}}E(e^{y},e^{Y})|^{2} dy  \ll
 \int_{\log 2}^{Y}  \left( \frac{y^{2}e^{y}}{e^{2Y}} +
      \frac{\frac{1}{y^{2}}e^{y}}{(e^{2Y})^{1 - \epsilon}} +
      \frac{Y}{e^{Y}} + \frac{1}{e^{y}} \right) dy \ll 1
\]
and we have
\begin{equation}
\begin{split}
   &  \int_{\log 2}^{Y} |\epsilon^{(T)}(y)|^{2} dy
     \ll  \int_{\log 2}^{Y} \left| \sum_{T \le \gamma \le e^{Y}}
    \frac{e^{i\gamma y}}{\rho \zeta^{'}(\rho)}
    \right|^{2} dy + O(1) \\
  & \le \sum_{j=0}^{\lfloor Y \rfloor}   \int_{\log 2 + j}^{\log 2 + j
  + 1}
    \left| \sum_{T \le \gamma \le e^{Y}}
    \frac{e^{i\gamma y}}{\rho \zeta^{'}(\rho)}
    \right|^{2} dy + O(1)
    \ll Y \frac{\log T}{T^{\frac{1}{4}}} + 1  \\
\end{split}
\end{equation}
where Lemma 6 has been applied in the last inequality.

By applying Lemmas 7-10, we may now prove Theorem 2.

\noindent {\it Proof of Theorem 2.} Once again the proof follows
Theorem 1.1 of \cite{RS} pp.\,180-181.
Let $f$ be a Lipschitz bounded continuous function that satisfies \\
$ |f(x) - f(y)| \le c_{f} |x - y | $.  By an application of the
Lipschitz condition, Cauchy-Schwarz, and Lemma 10, we have
\begin{equation}
\begin{split}
   \frac{1}{Y} \int_{\log 2}^{Y} f(\phi(y)) \, dy
   & = \frac{1}{Y} \int_{\log 2}^{Y} f(\phi^{(T)}(y)) \, dy
   + O \left( \frac{c_{f}}{Y}
   \int_{\log 2}^{Y} |\epsilon^{(T)}(y)| \, dy \right)  \\
   & = \frac{1}{Y} \int_{\log 2}^{Y} f(\phi^{(T)}(y)) \, dy
   + O \left( \frac{c_{f}}{\sqrt{Y}} \left(
   \int_{\log 2}^{Y} |\epsilon^{(T)}(y)|^{2} \, dy \right)^{\frac{1}{2}}
   \right) \\
  & = \frac{1}{Y} \int_{\log 2}^{Y} f(\phi^{(T)}(y)) \, dy
   + O \left( c_{f} \left( \frac{\log T}{T^{\frac{1}{4}}}
     + \frac{1}{\sqrt{Y}} \right)^{\frac{1}{2}}
     \right) . \\
\end{split}
\end{equation}
By Lemma 8, there is a distribution function $\nu_{T}$ for each
$T \ge \gamma_{1}$ such that
\[
   \nu_{T}(f) = \int_{\mathbb{R}}  f(x) \,  d\nu_{T}(x) =
   \lim_{Y \to \infty}  \frac{1}{Y} \int_{\log 2}^{Y} f(\phi^{(T)}(y))
   \, dy \ .
\]
Taking limits as $Y \to \infty$ we deduce that
\begin{equation}
\begin{split}
 & \nu_{T}(f)
 - O \left( \frac{c_{f}(\log T)^{\frac{1}{2}}}{T^{\frac{1}{8}}} \right)
  \le \liminf_{Y \to \infty} \frac{1}{Y} \int_{\log 2}^{Y} f(\phi(y)) dy \\
 & \le \limsup_{Y \to \infty} \frac{1}{Y} \int_{\log 2}^{Y} f(\phi(y)) dy
 \le \nu_{T}(f)
 + O \left( \frac{c_{f}(\log T)^{\frac{1}{2}}}{T^{\frac{1}{8}}} \right)
   \ . \\
   \label{eq:ublb}
\end{split}
\end{equation}
By Lemma 9$\, (i)$, we may choose a subsequence $\nu_{T_{k}}$
of these distribution functions $\nu_{T}$ and a
generalized distribution function $\nu$ such that
$\nu_{T_{k}} \to \nu$ weakly. By Lemma 9$\, (ii)$
\[
    \nu_{T_{k}}(f) = \int_{\mathbb{R}} f(x) \, d\nu_{T_{k}}(x)  \to
    \int_{\mathbb{R}} f(x) \, d\nu(x) = \nu(f)  \ .
\]
Replacing $T$ by $T_{k}$ and letting $k \to \infty$ in~(\ref{eq:ublb})
, we observe that
\[
   \lim_{Y \rightarrow \infty}
    \frac{1}{Y} \int_{\log 2}^{Y} f(\phi(y)) \, dy
   = \int_{\mathbb{R}} f(x) \, d\nu(x)  \ .
\]
Thus~(\ref{eq:ublb}) becomes
\begin{equation}
   \left|  \int_{\mathbb{R}} f(x) d\nu(x)
          - \int_{\mathbb{R}} f(x) d\nu_{T}(x) \right|
   \ll \frac{c_{f} (\log T)^{\frac{1}{2}}}{T^{\frac{1}{2}}} \ .
   \label{eq:ublb2}
\end{equation}
However, by applying equation~(\ref{eq:ublb2}) with $f(x)=1$
\[
   \left| \int_{\mathbb{R}} \ d \nu(x) - 1 \right|
   \ll  \frac{(\log T)^{\frac{1}{2}}}{T^{\frac{1}{8}}}
\]
and we conclude that $\nu$ is a distribution function by letting
$T \to \infty$.

By assuming the linear independence conjecture, we may provide a
concrete description of the Fourier transform of $\nu$ in terms of
the zeros of $\zeta(s)$.  This description will be practical in
obtaining finer details regarding $M(x)$.

\noindent {\it Proof of Corollary 1.}
The Fourier transform of $\nu$ is
\[ \widehat{\nu}(\xi) = \int_{\mathbb{R}} e^{-i \xi t} \ d \nu(t) \ .
\]
In the proof of Theorem 2, we demonstrated $\nu_{T} \to \nu$ weakly.
Hence, by Levy's Theorem (Lemma 9$\, (iii)$),
$\widehat{\nu}_{T} \to \widehat{\nu}$.
By Lemmas 7 and 8, we know that $\nu_{T}$ is
constructed from normalized Haar measure $\mu$ on the torus
$A \subset \mathbb{T}^{N}$ where $A$ is the topological closure of the
set
$\{ \,
   \left( \{ y \frac{\gamma_{1}}{2 \pi} \} ,
   \ldots , \{ y \frac{\gamma_{N}}{2 \pi} \}  \right)
 \ | \ y \in \mathbb{R} \ \} $.
However, the assumption of LI implies by the Kronecker-Weyl
theorem that $A = \mathbb{T}^{N}$ and  consequently normalized
Haar measure $d\mu = d\theta_{1} \ldots d\theta_{N}$ is Lesbesgue
measure on $\mathbb{T}^{N}$.  Hence, we observe
by~(\ref{eq:pullback}) and the change of variable formula for
integrals that $\widehat{\nu}_{T}(\xi)$ equals
\[
   \int_{\mathbb{R}} e^{-i \xi t} \, d\nu_{T}(t)
   = \int_{\mathbb{T}^{N}}
                        e^{-i\xi X_{T}(\theta)} d\mu
   = \int_{\mathbb{T}^{N}}
   e^{-i \xi \sum_{j=1}^{N} 2 \mathrm{Re} \left(
   \frac{1}{\rho \zeta^{'}(\rho)} e^{2 \pi i \theta_{j}}
   \right)} \ d\theta_{1} \ldots d\theta_{N} \\
\]
and it follows that
\[
  \widehat{\nu}(\xi) = \lim_{T \to \infty} \widehat{\nu}_{T} (\xi)
  = \lim_{T \to \infty} \prod_{j=1}^{N} \int_{0}^{1}
  e^{-i \xi 2 \mathrm{Re} \left(
   \frac{1}{\rho \zeta^{'}(\rho)} e^{2 \pi i \theta}
   \right)} \ d \theta \ .
\]
However the integral within the product equals
\[
  \int_{0}^{1}
  e^{-i \xi 2 \mathrm{Re} \left(
  \frac{1}{|\rho \zeta^{'}(\rho)|} e^{2 \pi i (\theta - \alpha_{\gamma})}
  \right)} \ d \theta
  = \int_{0}^{1}
  e^{-i \xi 2  \left(
  \frac{1}{|\rho \zeta^{'}(\rho)|} \cos{2 \pi  \theta}
  \right)} \ d \theta
\]
where $\alpha_{\gamma} = \mathrm{arg} (\rho \zeta^{'}(\rho))/2\pi$
and the last step follows by the periodicity of the integrand.
From the well-known identity for the $\tilde{J}_{0}$ Bessel function
\[ \int_{0}^{1} e^{i s \cos(2 \pi x)} dx =
   \frac{1}{\pi} \int_{0}^{\pi} \cos(s \sin x) dx
   = \tilde{J}_{0}(s)
\]
it follows that
\[
   \widehat{\nu}(\xi) = \prod_{\gamma > 0} \tilde{J}_{0}
   \left( \frac{2 \xi}{|(\frac{1}{2}+i\gamma)
   \zeta^{'}(\frac{1}{2}+i\gamma )|} \right) \ .
\]
We improve Theorem 1$(iv)$ by following closely Cram\'{e}r's
argument \cite{Cr}.

\noindent {\it Proof of Theorem 3}. Recall that By Lemma 5, we
have the decomposition
\begin{equation}
     M(e^{y})e^{-\frac{y}{2}}
      = \phi^{(T)}(y) + \epsilon^{(T)}(y)
     \label{eq:epf}
\end{equation}
where
\begin{equation}
  \phi^{(T)}(y) =
  \sum_{|\gamma| \le T} \frac{e^{iy\gamma}}{\rho \zeta^{'}(\rho)}
  \ , \
  \epsilon^{(T)}(y) = \sum_{T < |\gamma| \le e^{Y}}
     \frac{e^{iy\gamma}}{\rho \zeta^{'}(\rho)}
       + e^{-\frac{y}{2}} E(e^{y},e^{Y})
\end{equation}
and $E(x,T)$ is defined in~(\ref{eq:erm}). Consequently, we deduce
\begin{equation}
\begin{split}
  m(Y) & := \frac{1}{Y} \int_{0}^{Y} \left( \frac{M(e^{y})}{e^{\frac{y}{2}}}
  \right)^{2} \, dy
  = \frac{1}{Y} \int_{0}^{Y} |\phi^{(T)}(y)|^{2}
  \, dy + \frac{1}{Y} \int_{0}^{Y} |\epsilon^{(T)}(y)|^{2}
  \, dy  \\
  & + O \left( \left( \frac{1}{Y} \int_{0}^{Y} |\phi^{(T)}(y)|^{2}
  \right)^{\frac{1}{2}}
  \left( \frac{1}{Y} \int_{0}^{Y} |\epsilon^{(T)}(y)|^{2}
  \right)^{\frac{1}{2}}
  \right) \ . \\
  \label{eq:mean}
\end{split}
\end{equation}
As the second integral was treated in Lemma 10,
we concentrate on the first integral in~(\ref{eq:mean}).
Squaring out the terms in $\phi^{(T)}(y)$, we deduce
\begin{equation}
\begin{split}
   \int_{1}^{Y} |\phi^{(T)}(y)|^{2} \, dy
   & = (Y-1)\sum_{\gamma \le T} \frac{2}{|\rho \zeta^{'}(\rho)|^{2}}
   \\
   & +   \sum_{{\begin{substack}{0 < |\gamma|, |\lambda| < T
         \\ \gamma \ne \lambda  }\end{substack}}}
   \frac{1}{(\frac{1}{2}+i\gamma)
    \zeta^{'}(\rho)(\frac{1}{2} +i \lambda) \zeta^{'}(\rho^{'})}
   \int_{1}^{Y} e^{iy(\gamma +\lambda)} \, dy  \ .  \\
\end{split}
\end{equation}
In the second sum, the contribution from pairs $(\gamma, \lambda)$
with the same sign is
\[
  \sum_{0 < \gamma, \lambda \le T}
  \frac{1}{\gamma |\zeta^{'}(\rho)| \lambda |\zeta^{'}(\rho^{'})|
  (\gamma + \lambda)} \ll
  \left( \sum_{0 < \gamma < T}
  \frac{1}{\gamma^{\frac{3}{2}}|\zeta^{'}(1/2+i\gamma)|}
  \right)^{2} \ll 1  \ .
\]
Here we have applied  $x+y \ge 2\sqrt{xy}$ and then evaluated
the resulting sum by a partial summation similar to Lemma 1$(iii)$.
Also note that
\[
  \sum_{\gamma < T} \frac{1}{|\rho \zeta^{'}(\rho)|^{2}}
  = \beta - \sum_{\gamma > T}  \frac{1}{|\rho \zeta^{'}(\rho)|^{2}}
  = \beta + O \left( \frac{1}{T} \right)
\]
where $\beta$ is defined by~(\ref{eq:zsum}) and the error term is obtained
by Lemma 1$(ii)$.
We have now shown that
\[
   \frac{1}{Y}  \int_{1}^{Y} |\phi^{(T)}(y)|^{2} \, dy
   = \beta + O \left( \frac{1}{T} + \frac{1}{Y}  +
   \frac{\Sigma(T)}{Y} \right)
\]
where
\begin{equation}
\begin{split}
   \Sigma = \Sigma(T,Y) & =
         \sum_{{\begin{substack}{0 < \gamma, \lambda < T
         \\ \gamma \ne \lambda  }\end{substack}}}
   \frac{1}{\gamma
   |\zeta^{'}(\rho)|\lambda
   |\zeta^{'}(\rho^{'})|}
   \min \left( Y, \frac{1}{|\gamma-\lambda|} \right) \, dy  \\
   & = \Sigma_{1}(T,Y) + \Sigma_{2}(T,Y) \ . \\
\end{split}
\end{equation}
The first sum is the contribution from those pairs with
$|\gamma - \lambda| \le 1$ and the second sum consists of the
complementary terms.
We have
\begin{equation}
\begin{split}
  \Sigma_{2}(T,Y)
  & \le \sum_{0 < \gamma < T} \frac{1}{\gamma |\zeta^{'}(\rho)|}
      \left( \sum_{\lambda < \gamma^{\frac{1}{2}}}
           + \sum_{\gamma^{\frac{1}{2}} < \lambda < \gamma-\gamma^{\frac{1}{2}}}
           + \sum_{\gamma-\gamma^{\frac{1}{2}} < \lambda < \gamma-1} \right. \\
    & \hspace{3cm}  \left.
    + \sum_{\gamma+1 < \lambda < \gamma+ \gamma^{\frac{1}{2}}}
           + \sum_{\gamma+\gamma^{\frac{1}{2}} < \lambda < 2\gamma}
           + \sum_{2\gamma < \lambda}  \right)
    \frac{1}{\lambda |\zeta^{'}(\rho^{'})| |\gamma-\lambda|}  \ . \\
    & = \sigma_{1} +  \sigma_{2} +  \sigma_{3} +  \sigma_{4} +
  \sigma_{5} +  \sigma_{6} \\
\end{split}
\end{equation}
By a calculation completely analogous to the calculation in
Lemma 6, we obtain
\[
  \sigma_{1} \le \sum_{0 < \gamma < T}
  \frac{1}{\gamma(\gamma-\gamma^{\frac{1}{2}}) |\zeta^{'}(\rho)|}
  \ll
  \sum_{\gamma > 0} \frac{1}{\gamma^{2} |\zeta^{'}(\rho)|} \ll
  1 \ ,
\]
\[
  \sigma_{2} \le \sum_{0 < \gamma < T}
  \frac{1}{\gamma^{\frac{3}{2}} |\zeta^{'}(\rho)|}
  \left( \sum_{\lambda \ge \gamma^{\frac{1}{2}}} \frac{1}{
  |\lambda \zeta^{'}(\rho^{'})|^{2}} \right)^{\frac{1}{2}}
  (\gamma \log \gamma)^{\frac{1}{2}}
  \ll \sum_{\gamma > 0}
  \frac{(\log \gamma)^{\frac{1}{2}}}
  {\gamma^{\frac{5}{4}} |\zeta^{'}(\rho)|}
  \ll 1  \ ,
\]
\[
  \sigma_{3} \le \sum_{0 < \gamma < T}
  \frac{1}{\gamma |\zeta^{'}(\rho)|}
  \left( \sum_{ \gamma - \gamma^{\frac{1}{2}} < \lambda}
  \frac{1}{|\rho^{'} \zeta^{'}(\rho^{'})|^{2}} \right)^{\frac{1}{2}}
  (\gamma^{\frac{1}{2}} \log \gamma)^{\frac{1}{2}} \ll
  \sum_{\gamma > 0}
  \frac{(\log \gamma)^{\frac{1}{2}}}
  {\gamma^{\frac{5}{4}} |\zeta^{'}(\rho)|}
  \ll 1
\]
where we have applied Lemma 1$(iii)$ in each of these cases.
The computation of $\sigma_{4}$ is analogous to $\sigma_{3}$ and the
computation of $\sigma_{5}$ is analogous to $\sigma_{2}$
\[
  \sigma_{4} \ll \sum_{\gamma > 0}
  \frac{(\log \gamma)^{\frac{1}{2}}}
  {\gamma^{\frac{5}{4}} |\zeta^{'}(\rho)|}
  \ll 1     \ , \
   \sigma_{5} \ll \sum_{\gamma > 0}
  \frac{(\log \gamma)^{\frac{1}{2}}}
  {\gamma^{\frac{5}{4}} |\zeta^{'}(\rho)|}
  \ll 1     \ .
\]
For the final sum we obtain
\begin{equation}
\begin{split}
  & \sigma_{6}  \le
  \sum_{0 < \gamma < T} \frac{1}{\gamma |\zeta^{'}(\rho)|}
  \sum_{k\ge 1} \frac{1}{(2^{k}-1)\gamma}
    \sum_{2^{k}\gamma \le \lambda \le 2^{k+1}\gamma}
  \frac{1}{\lambda |\zeta^{'}(\frac{1}{2}+i \lambda)|} \\
  &
  \ll
   \sum_{0 < \gamma < T} \frac{1}{\gamma^{2} |\zeta^{'}(\rho)|}
  \sum_{k \ge 1} \frac{1}{2^{k}}
  \left( \sum_{2^{k}\gamma  < \lambda < 2^{k+1} \gamma}
  \frac{1}{|\lambda \zeta^{'}(\rho^{'})|^{2}}
  \right)^{\frac{1}{2}} ((2^{k}\gamma) \log(2^{k}
  \gamma))^{\frac{1}{2}} \\
  & \ll
   \sum_{0 < \gamma < T} \frac{1}{\gamma^{2} |\zeta^{'}(\rho)|}
  \sum_{k \ge 1}
  \frac{(\log(2^{k} \gamma))^{\frac{1}{2}}}{2^{k}}
  \ll   \sum_{\gamma > 0}
  \frac{\sqrt{\log \gamma}}{\gamma^{2} |\zeta^{'}(\rho)|}
  \ll 1 \\
\end{split}
\end{equation}
and we deduce that
\[
  \Sigma_{2}(T,Y) \ll \sum_{\gamma > 0}
  \frac{1}{\gamma^{\frac{5}{4}}|\zeta^{'}(\rho)|}
  \ll 1 \ .
\]
Thus
\begin{equation}
   \frac{1}{Y}  \int_{1}^{Y} |\phi^{(T)}(y)|^{2} \, dy
   = \beta + O \left( \frac{1}{T} + \frac{1}{Y}  +
   \frac{\Sigma_{1}}{Y} \right)
   \label{eq:newb}
\end{equation}
where
\[
   \Sigma_{1} =
   \Sigma_{1}(T,Y) =
        \sum_{{\begin{substack}{0 < \gamma, \lambda < T
         \\ |\gamma-\lambda| \le 1  }\end{substack}}}
   \frac{1}{\gamma
   |\zeta^{'}(\rho)|\lambda
   |\zeta^{'}(\rho^{'})|}
   \min \left( Y, \frac{1}{|\gamma-\lambda|} \right)  \ .
\]
In addition, we know by Lemma 10 that
\begin{equation}
  \frac{1}{Y} \int_{1}^{Y} |\epsilon^{(T)}(y)|^{2} \, dy
  \ll \frac{\log T}{T^{\frac{1}{4}}} + \frac{1}{Y} \ .
  \label{eq:le7}
\end{equation}
Let $0 < \eta < 1$. Choose and fix $T=T_{\eta}$ large enough to
make the $O(\frac{1}{T})$  in~(\ref{eq:newb}) and $O(\frac{\log
T}{T^{\frac{1}{4}}})$ in~(\ref{eq:le7}) less than $\eta$.  Choose
$Y_{1}$ large enough such that if $Y\ge Y_{1}$ the $O(Y^{-1})$
expressions in~(\ref{eq:newb}) and~(\ref{eq:le7}) are less than
$\eta$. Choose $Y_{\eta}$ to satisfy
\begin{equation}
   Y_{\eta} = \max \left( \frac{1}{\eta \min_{0<\gamma \le
         T_{\eta}}
   |\gamma^{'} -
   \gamma|} , Y_{1} \right)
\end{equation}
where if $\gamma$ denotes an imaginary ordinate of a zero of
$\zeta(s)$ then $\gamma^{'}$ is the next largest one
(note that $\gamma^{'} > \gamma$ since $J_{-1}(T) \ll T$
implies all zeros are simple).
We will consider $Y \ge Y_{\eta}$ and analyze $\Sigma_{1}$.
Decompose $\Sigma_{1}(T_{\eta},Y) =
\Sigma_{11}(T_{\eta},Y)+\Sigma_{12}(T_{\eta},Y)$
where  the first sum contains pairs $(\gamma,\lambda)$ with
$|\gamma-\lambda|^{-1}
\le \eta Y$ and the second sum contains the complementary set.
Therefore
\begin{equation}
\begin{split}
  & \Sigma_{11}(T_{\eta},Y)
   \le \eta Y \sum_{\gamma < T_{\eta}} \frac{1}{\gamma
  |\zeta^{'}(\rho)|}
  \sum_{\gamma-1 < \lambda < \gamma+1} \frac{1}{\lambda
    |\zeta^{'}(\rho^{'})|} \\
  & \ll  \eta Y \sum_{\gamma < T_{\eta}} \frac{1}{\gamma
  |\zeta^{'}(\rho)|}  \left(
  \sum_{\lambda > \gamma-1} \frac{1}{|\lambda
  \zeta^{'}(\rho^{'})|^{2}}  \right)^{\frac{1}{2}}
  (\log \gamma)^{\frac{1}{2}}
   \ll
  \eta Y \sum_{\gamma < T_{\eta}}
  \frac{(\log \gamma)^{\frac{1}{2}}}{\gamma^{\frac{3}{2}}
  |\zeta^{'}(\rho)|}
\end{split}
\end{equation}
and we have $ \Sigma_{11}(T_{\eta},Y) \le c_{3} \eta Y $ for
$c_{3} > 0$ by Lemma 1$(iii)$. The second sum consists of pairs
$(\gamma,\lambda)$ such that
\[
   \eta Y < |\gamma - \lambda|^{-1}
   < \left( \min_{0<\gamma \le
         T_{\eta}} |\gamma^{'} - \gamma| \right)^{-1}
\]
which implies
\[
    Y < \left( \eta \min_{0<\gamma \le
         T_{\eta}} |\gamma^{'} - \gamma| \right)^{-1}
      \le Y_{\eta}
\]
and thus this second sum is empty.
Consequently, $\Sigma_{12}(T_{\eta},Y)
 = 0$ and thus
$ \Sigma_{1}(T_{\eta},Y)Y^{-1}  \le c_{3} \eta$.
This demonstrates that
\begin{equation}
   \left| \frac{1}{Y}  \int_{1}^{Y} |\phi^{(T_{\eta})}(y)|^{2} \, dy
   - \beta \right| \le (2+c_{3}) \eta \ \ and \ \
     \frac{1}{Y}  \int_{1}^{Y} |\epsilon^{(T_{\eta})}(y)|^{2} \, dy
   \le  2 \eta \
   \label{eq:bounds}
\end{equation}
if $Y \ge Y_{\eta}$.
By~(\ref{eq:mean}) and~(\ref{eq:bounds}) we deduce
\[
   \left| \frac{1}{Y} \int_{1}^{Y} \left(
   \frac{M(e^{y})}{e^{\frac{y}{2}}}  \right)^{2} \, dy - \beta \right|
   \le (4+c_{3}) \eta + c_{4} \sqrt{\eta}
   \ll \sqrt{\eta}
\]
if $Y \ge Y_{\eta}$ and hence the proof is finished.
\section{Applications of LI}
The goal of this section is to study the true order of $M(x)$.
We will attempt to find the size of the tail of the probability
measure $\nu$ associated to $\phi(y) = e^{-\frac{y}{2}}M(e^{y})$.
The tool we employ in studying tails of $\nu$ are
probability results due to Montgomery \cite{Mo}.
We will need
to assume the linear independence conjecture for our analysis.
Consider a random variable $X$, defined on the infinite torus
$\mathbb{T}^{\infty}$ by
\[ X(\underline{\theta}) = \sum_{k = 1}^{\infty} r_{k} \sin(2 \pi \theta_{k})
\]
where $\underline{\theta} = (\theta_{1},\theta_{2}, \ldots ) \in
\mathbb{T}^{\infty}$
and $r_{k} \in \mathbb{R}$ for $k \ge 1$.  This
 is a map $X : \mathbb{T}^{\infty} \rightarrow \mathbb{R} \cup \{
 \infty \}$.
Under the assumption
$\sum_{k \ge 1} r_{k}^{2} < \infty$, Komolgorov's theorem ensures that
$X$ converges almost everywhere.
In addition, $\mathbb{T}^{\infty}$ possesses a canonical probability measure
$P$.  Attached to the random variable $X$ is the distribution function
$\nu_{X}$ defined by
\[ \nu_{X}(x) = P(X^{-1}(-\infty,x)) .
\]
For these random variables, Montgomery \cite{Mo} pp.\,14-16 proved
the following results.
\newtheorem{mont}[asym]{Lemma}
\begin{mont}
Let $X(\underline{\theta}) = \sum_{k=1}^{\infty} r_{k} \sin 2
\pi \theta_{k}$ where
$\sum_{k=1}^{\infty} r_{k}^{2} < \infty$.
For any integer
$K \ge 1$, \\
(i)
\[ P\left(X(\underline{\theta})  \ge 2 \sum_{k=1}^{K} r_{k} \right)
   \le \exp \left( -\frac{3}{4} \left(\sum_{k=1}^{K} r_{k}\right)^{2}
    \left(\sum_{k>K} r_{k}^{2}\right)^{-1} \right)
\]
(ii) and if $\delta$ is so small that $\sum_{r_{k} > \delta} (r_{k} - \delta)
\ge V$, then
\[ P(X(\underline{\theta}) \ge V) \ge \frac{1}{2}
  \exp \left(-\frac{1}{2} \sum_{r_{k} > \delta}
  \log \left( \frac{\pi^{2} r_{k}}{2 \delta} \right) \right) \ .
\]
\end{mont}
Observe that the linear independence assumption implies that the limiting
distribution $\nu$ constructed in Theorem 2 equals
$\nu_{X}$ where $X$ is the random variable
\[ X(\underline{\theta})
   = \sum_{\gamma > 0}
   \frac{2}{|\rho \zeta^{'}(\rho)|} \sin(2 \pi \theta_{\gamma}) \ .
\]
In the above sum $\gamma$ ranges over the positive imaginary ordinates of the
zeros of $\zeta(s)$.  We abbreviate notation by setting
$r_{\gamma} = \frac{2}{|\rho \zeta^{'}(\rho)|}$.
By assuming the linear independence conjecture, we
may now study $\nu$ via the random variable $X$.
By applying Lemma 11,
we can estimate the tails of the limiting distribution $\nu$.
Define
\[
   a(T) := \sum_{\gamma < T} r_{\gamma}
        = \sum_{\gamma < T} \frac{2}{|\rho \zeta^{'}(\rho)|}
        \ \ and \ \
   b(T) := \sum_{\gamma \ge T} r_{\gamma}^{2}
        = \sum_{\gamma \ge T} \frac{4}{|\rho \zeta^{'}(\rho)|^{2}}   \ .
\]
By Lemma 1, the conjectured
formulae are
\begin{equation}
    a(T) \asymp (\log T)^{\frac{5}{4}} \
   \ and \ \ b(T)
   \asymp \frac{1}{T} \ .
   \label{eq:quitte}
\end{equation}
Assuming these bounds we prove upper and lower bounds for
the tail of the limiting distribution $\nu$.
Let $V$ be a large parameter. Our goal is to find upper and
lower bounds for the tail of the probability distribution, namely
\[
   \nu([V,\infty)) := \int_{V}^{\infty} d\nu(x) =
   P( X(\underline{\theta}) \ge V) \ .
\]
\subsection{An upper bound for the tail}
Choose $T$ such that $a(T^{-}) < V \le a(T)$.
Note that $T$ is the ordinate of a zero.  We have the chain of inequalities
\begin{equation}
  (\log T)^{\frac{5}{4}} \ll a(T^{-}) < V
  \le a(T) \ll (\log T)^{\frac{5}{4}}.
  \label{eq:fuir}
\end{equation}
This implies $\log T \asymp V^{\frac{4}{5}}$ and we have by Lemma
11$(i)$,~(\ref{eq:quitte}), and~(\ref{eq:fuir}),
\begin{equation}
\begin{split}
    P \left( X(\underline{\theta}) \ge c_{5}V \right)
   & \le
   P \left( X(\underline{\theta}) \ge 2a(T) \right)  \le
   \exp \left( -\frac{3}{4}a(T)^{2}b(T)^{-1} \right) \\
   & \le \exp \left(-c_{6} V^{2} T \right)
   \le \exp \left(-c_{6} V^{2}
   e^{ \left( c_{7} V \right)^{\frac{4}{5}}} \right)  \\
\end{split}
\end{equation}
for effective constants $c_{5},c_{6},$ and $c_{7}$. By altering
the constants, we obtain the upper bound
\[
   P( X(\underline{\theta}) \ge V)
   \ll  \exp (- \exp(c_{7} V^{\frac{4}{5}})) \ .
\]
\subsection{A lower bound for the tail}
This is a more delicate analysis than the upper bound.
We now apply Lemma 11$(ii)$.  As before, $V$ is
considered fixed and large. We would like to choose
$\delta$ small enough such that
\begin{equation}
   \sum_{r_{\gamma} > \delta} (r_{\gamma} - \delta)
   \ge V \ .
   \label{eq:lower}
\end{equation}
Introduce the notation $S_{\delta}$ and $N_{\delta}$ such that
\[ S_{\delta} =  \{ \gamma \ | \ r_{\gamma} > \delta \ \} \
   and \
   N_{\delta} = \#  S_{\delta}
\]
where $\gamma$ ranges over positive imaginary ordinates of zeros
of $\zeta(s)$.
Let $\epsilon$ be a small fixed number.
Note that RH implies $|\zeta^{'}(\rho)| \ll
 |\rho|^{\epsilon}$. Thus,
\[ \delta < \frac{2}{c_{8}|\rho|^{1+\epsilon}}
   \Longrightarrow \delta <  \frac{2}{|\rho \zeta^{'}(\rho)|} \
\]
for some effective constant $c_{8}$.
However, notice that
\[ \delta < \frac{2}{c_{8}|\rho|^{1+\epsilon}}
   \Longleftrightarrow
   |\rho| \le \left( \frac{2}{c_{8} \delta}
   \right)^{\frac{1}{1+\epsilon}}
\]
and since $|\rho| \ll \gamma$, we obtain
\[
  \gamma \le c_{9}
  \left( \frac{1}{\delta} \right)^{\frac{1}{1+\epsilon}}
  \Longrightarrow \delta <  \frac{2}{|\rho \zeta^{'}(\rho)|} .
\]
We deduce from Riemann's zero counting formula that there
are at least
\[
          c_{9} \left( \frac{1}{\delta} \right)^{\frac{1}{1+\epsilon}}
          \log \left( \frac{1}{\delta} \right)
          + O \left( \left( \frac{1}{\delta} \right)^{\frac{1}{1+\epsilon}}
            \right)
\]
zeros in the set $S_{\delta}$.
We will now find an upper bound for $N_{\delta}$. Gonek \cite{G2} has
defined the number
\[ \Theta = \mathrm{l.u.b.}
   \{ \ \theta \ | \ |\zeta^{'}(\rho)|^{-1}
      \ll |\gamma|^{\theta}, \ \forall \rho \ \}  .
\]
However $J_{-1}(T) \ll T$
implies $\Theta \le \frac{1}{2}$.
Gonek \cite{G2} has speculated that
$\Theta = \frac{1}{3}$.  Choose $\epsilon < \frac{1}{2}$.
This implies that if $\gamma \in S_{\delta}$
then
\[ \delta < \frac{2}{|\rho \zeta^{'}(\rho)|}
   \ll \frac{|\rho|^{\frac{1}{2} + \epsilon}}{|\rho|}
   \le \frac{1}{|\gamma|^{\frac{1}{2} - \epsilon}} \ .
\]
We deduce that if $\gamma \in S_{\delta}$
then $\gamma \ll \left( \frac{1}{\delta} \right)^{2 + \epsilon}$
where $\epsilon$ has been taken smaller.
We conclude that $N_{1}(\delta) \le N_{\delta} \le N_{2}(\delta)$ where
\[
   N_{1}(\delta) =
   c_{9} \left( \frac{1}{\delta} \right)^{1-\epsilon}
   \ and \
  N_{2}(\delta) =   c_{10} \left(
   \frac{1}{\delta} \right)^{2 + \epsilon}  \ .
\]
We are trying to determine a condition on $\delta$
so that~(\ref{eq:lower}) will be satisfied.
Note that
\[ \sum_{r_{\gamma} > \delta} (r_{\gamma} - \delta)
   \ge
   \sum_{\gamma \le N_{1}} (r_{\gamma} - \delta) \ .
\]
Before evaluating the second sum, observe that
\[ \delta N_{1} = c_{9}
   \delta^{\epsilon}
     \rightarrow 0
   \ as \
   \delta \rightarrow 0  .
\]
We will choose $\delta$ as a function of $V$ and as $V \rightarrow
\infty$ we have $\delta \rightarrow 0$. However,
by~(\ref{eq:quitte})
\begin{equation}
\begin{split}
   & \sum_{\gamma \le N_{1}} (r_{\gamma} - \delta)
    = 2 \sum_{\gamma \le N_{1}} \frac{1}{|\rho \zeta^{'}(\rho)|}
   - \delta \sum_{\gamma \le N_{1}} 1 \\
   & \ge c_{11} (\log N_{1})^{\frac{5}{4}}
     - \frac{\delta N_{1}}{2 \pi} \log N_{1} + O(\delta N_{1})
   \ge c_{12} (\log N_{1})^{\frac{5}{4}} \\
\end{split}
\end{equation}
where $0 < c_{12} < c_{11}$.
The last inequality holds for $N_{1}$ sufficiently large.
Hence, choosing $N_{1} =
\exp(( V/c_{12})^{\frac{4}{5}})$ implies
\[ \sum_{r_{\gamma} > \delta} (r_{\gamma} - \delta)
   \ge   \sum_{\gamma \le N_{1}} (r_{\gamma} - \delta)
   \ge V \ .
\]
Thus if $\delta$ satisfies
\[ c_{9} \left( \frac{1}{\delta} \right)^{1-\epsilon}
   = \exp(( V/c_{12})^{\frac{4}{5}})
\]
(i.e. $\delta = c_{13} \exp( -c_{14} V^{\frac{4}{5}})$) we have
satisfied~(\ref{eq:lower}). By this choice of $\delta$, Lemma
11$(ii)$ implies
\begin{equation}
  P(X(\underline{\theta}) \ge V) \ge \frac{1}{2}
  \exp \left(-\frac{1}{2} \sum_{r_{\gamma} > \delta}
  \log \left( \frac{\pi^{2} r_{\gamma}}{2 \delta} \right) \right) \ .
  \label{eq:lem2}
\end{equation}
An upper bound of the sum will provide a lower bound for the tail.
Note that $\frac{1}{|\rho \zeta^{'}(\rho)|} \rightarrow 0$ under
the assumption that all zeros are simple (see \cite{T}
pp.\,377-380). Consequently $\frac{1}{|\rho \zeta^{'}(\rho)|} \le
c_{15}$ and we obtain
\[
   \sum_{r_{\gamma} > \delta}
   \log \left( \frac{\pi^{2} r_{\gamma}}{2 \delta} \right)
   \le
   \sum_{\gamma \le N_{2}}
   \log \left( \frac{\pi^{2} c_{15}}{\delta} \right)
    \ll
   \log \left( \frac{\pi^{2} c_{15}}{\delta} \right)
   N_{2} \log N_{2} \ . \\
\]
By definition of $N_{2}(\delta)$ and our choice of $\delta$ it
follows that
\begin{equation} \sum_{r_{\gamma} > \delta}
   \log \left( \frac{\pi^{2} r_{\gamma}}{2 \delta} \right)
   \ll V^{\frac{4}{5}} \exp( c_{16}V^{\frac{4}{5}} )
   \ll \exp( c_{17}V^{\frac{4}{5}} ) \ .
   \label{eq:rglb}
\end{equation}
By~(\ref{eq:lem2}) and~(\ref{eq:rglb}) we arrive at the lower
bound
\[ P(X(\underline{\theta}) \ge V)
   \gg \exp \left( -
   \exp( c_{18}V^{\frac{4}{5}} )
   \right) \ .
\]
Putting this all together establishes the following
highly conditional result.
\newtheorem{tai}[asym]{Corollary}
\begin{tai}
The Riemann hypothesis, the linear independence conjecture,
\[
   \sum_{0 < \gamma < T} \frac{1}{|\rho \zeta^{'}(\rho)|}
   \asymp (\log T)^{\frac{5}{4}} \ , \ and \
   \sum_{\gamma > T} \frac{1}{|\rho \zeta^{'}(\rho)|^{2}}
   \asymp \frac{1}{T}
\]
imply
\[
   \exp(-
   \exp(\tilde{c}_{1}V^{\frac{4}{5}})) \ll
           \nu([V,\infty)) \ll
   \exp(-
   \exp(\tilde{c}_{2}V^{\frac{4}{5}}))
\]
for effective constants $\tilde{c}_{i} > 0$ for $i=1 \ldots 2$.
\end{tai}
\subsection{Speculations on the lower order of $M(x)$}
We now examine the effect that
bounds for the tail of the probability measure
have on the lower order of $M(x)$.
Note that the following argument is only heuristic.
We begin with the lower bound
\[ \exp \left( -
   \exp( \tilde{c}_{1} V^{\frac{4}{5}} ) \right)
   \ll \nu([V,\infty)) \ .
\]
Assuming the linear independence conjecture,
the Riemann hypothesis, and $J_{-1}(T) \ll T$,  we have
\begin{equation}
  \lim_{Y \rightarrow \infty}
   \frac{1}{Y} \mathrm{meas}
   \{ y \in [0,Y] \ | \ M(e^{y}) \ge e^{y/2}V \ \} =
   \nu([V,\infty)) \ .
   \label{eq:dog}
\end{equation}
We will assume that the convergence of~(\ref{eq:dog})
is sufficiently uniform in $Y$.
By~(\ref{eq:dog}) there exists a function $f(V)$,  such that
\[
    \frac{1}{Y} \mathrm{meas}
   \{ y \in [0,Y] \ | \ M(e^{y}) \ge e^{y/2}V \ \}  \gg\exp \left( -
   \exp( \tilde{c}_{1}V^{\frac{4}{5}} )
   \right)
\]
if $V$ is sufficiently large and $Y \ge f(V)$.  We now choose $Y$ as
a function of $V$ by the equation
\[
    V =  \left(
   \frac{\theta}{\tilde{c}_{1}} \right)^{\frac{5}{4}}
   ( \log_{2} Y )^{\frac{5}{4}} \ \mathrm{or} \
   Y = g(V) = \exp \left( \exp \left( \frac{\tilde{c}_{1}}{\theta}
   V^{\frac{4}{5}} \right)
   \right)
\]
for $0 < \theta < 1$. If we had $g(V) \ge f(V)$ then it
would follow that for large $Y$,
\[
  \exp \left( \log Y - ( \log Y)^{\theta}
   \right)
   \ll
   \mathrm{meas}
   \{ y \in [0,Y] \ | \ M(e^{y})e^{-\frac{y}{2}} \ge \alpha
   (\log_{2} Y )^{\frac{5}{4}} \ \}
\]
where $\alpha = \left( \frac{\theta}{\tilde{c}_{1}} \right)^{\frac{5}{4}}$.
Since $0 < \theta < 1$
 the left hand side of the equation approaches infinity
as $Y \to \infty$.
In turn, this implies that there exists an increasing sequence
of real numbers $y_{m}$ such that $y_{m} \rightarrow \infty$
and
\[ \frac{M(e^{y_{m}})}{e^{\frac{y_{m}}{2}}} \ge
   \alpha \, ( \log_{2} y_{m} )^{\frac{5}{4}} .
\]
Suppose by way of contradiction, that the above inequality is
false.  That is, there exists a real number $u_{0}$ such that
\[ \frac{M(e^{y})}{e^{\frac{y}{2}}}
   < \alpha \,
   ( \log_{2} y )^{\frac{5}{4}}
\]
for all $y \ge u_{0}$.   Then we have that
\begin{equation}
\begin{split}
  &  \mathrm{meas}
   \{ y \in [0,Y] \ | \ M(e^{y}) \ge \alpha \, e^{\frac{y}{2}}
   ( \log_{2} Y )^{\frac{5}{4}} \ \} \\
   = \ &  \mathrm{meas}
   \{ y \in [0,u_{0}] \ | \ M(e^{y}) \ge \alpha \, e^{\frac{y}{2}}
   ( \log_{2} Y )^{\frac{5}{4}} \ \}  \\
\end{split}
\end{equation}
since if $u_{0} \le y \le Y$ then
\[
   \frac{M(e^{y})}{e^{\frac{y}{2}}}
   \le \alpha \,
   ( \log_{2} y )^{\frac{5}{4}}
   \le \alpha \,
   ( \log_{2} Y )^{\frac{5}{4}} .
\]
Thus we deduce that
\begin{equation}
   \exp(\log Y - (\log Y)^{\theta}) \le u_{0} \ll 1
\end{equation}
which is a contradiction for large enough $Y$.
Hence, our original assumption is false and we obtain
\[
   \limsup_{y \to \infty}
   \frac{M(e^{y})}{e^{\frac{y}{2}}(\log \log y)^{\frac{5}{4}}}
   \ge \left( \frac{\theta}{\tilde{c}_{1}} \right)^{\frac{5}{4}} \ .
\]
Letting $\theta \to 1$ we have
\[
    \limsup_{y \to \infty}
   \frac{M(e^{y})}{e^{\frac{y}{2}}(\log \log y)^{\frac{5}{4}}}
   \ge \left( \frac{1}{\tilde{c}_{1}} \right)^{\frac{5}{4}} \ .
\]
We now consider the upper bound.
Arguing in the same fashion we  have
\begin{equation}
   \nu([V,\infty)) = P( \underline{\theta} \in \mathbb{T}^{\infty} \ | \
                         X(\underline{\theta}) \ge V )
   \ll  \exp \left( -
   \exp( \tilde{c}_{2}V^{\frac{4}{5}} )
   \right)  \ .
   \label{eq:ldub}
\end{equation}
For $n \in \mathbb{N}$ define the event
\[
  A_{n} = \left\{ \underline{\theta} \in \mathbb{T}^{\infty} \ | \
             X(\underline{\theta}) \ge \left( \frac{1}{\tilde{c}_{2}}
             \log \log (n (\log n)^{\theta}) \right)^{\frac{5}{4}}
   \right\}
\]
with $\theta > 1$.
Therefore we have by~(\ref{eq:ldub})
\[
   \sum_{n=n_{0}}^{\infty} P(A_{n})
   \ll  \sum_{n=n_{0}}^{\infty} \frac{1}{n(\log n)^{\theta}}
   \ll 1
\]
for $n_{0}$ a sufficiently large integer.
By the Borel-Cantelli lemma, it follows that
\begin{equation}
  P(A_{n} \ infinitely \ often) = 0
\end{equation}
which suggests that if the convergence of~(\ref{eq:dog}) is
sufficiently uniform then
\[
        \limsup_{y \to \infty} \frac{M(e^{y})}{e^{\frac{y}{2}} (\log \log
  y)^{\frac{5}{4}}}
  \le \left( \frac{1}{\tilde{c}_{2}} \right)^{\frac{5}{4}} \ .
\]
Hence, our analysis shows that the bounds
\[
    \exp \left( - \exp( \tilde{c}_{1}V^{\frac{4}{5}} )
   \right) \ll
   \nu([V,\infty)) \ll
   \exp \left( -  \exp( \tilde{c}_{2}V^{\frac{4}{5}} )
   \right)
\]
suggest
\[
    \left( \frac{1}{\tilde{c}_{1}} \right)^{\frac{5}{4}} \le
        \limsup_{y \to \infty} \frac{M(e^{y})}{e^{\frac{y}{2}} (\log \log
  y)^{\frac{5}{4}}}
  \le \left( \frac{1}{\tilde{c}_{2}} \right)^{\frac{5}{4}} \ .
\]
Thus we arrive at an argument for the conjecture~(\ref{eq:lb}).

By the preceding heuristic analysis and Theorems 1-3
we hope to have demonstrated that the size of $M(x)$ depends in a crucial way
on the sizes of the discrete moments $J_{-\frac{1}{2}}(T)$ and $J_{-1}(T)$.

\noindent D\'epartement de Math\'ematiques et de statistique,
Universit\'e de Montr\'eal,
CP 6128 succ Centre-Ville,
Montr\'eal, QC, Canada  H3C 3J7 \\
EMAIL: nathanng@dms.umontreal.ca


\begin{thebibliography}{99}
%
%
%
\bibitem{B}
P.T. Bateman et al., {\em Linear relations connecting the imaginary
parts of the zeros of the zeta function},
in {\em Computers in Number Theory}, eds. A.O.L. Atkin and B.J. Birch,
Academic Press, New York, 1971, 11-19.

\bibitem{Bi}
Billingsley, P., {\em Probability and measure}, John Wiley and Sons,
New York, 1986.

\bibitem{C}
J.B. Conrey, {\em More than two-fifths of the zeros of the Riemann
zeta function are on the critical line}, J. Reine. Angew. Math.
\textbf{399} (1989), 1-26.

\bibitem{Cr}
Harald Cram\'{e}r, {\em Ein Mittelwertsatz in der Primzahltheorie},
Math. Z. \textbf{12} (1922), 147-153.

\bibitem{D}
Harold Davenport, {\em Multiplicative number theory, Second Edition},
Graduate Texts in Mathematics, Springer-Verlag, New York, 1980.

\bibitem{Ga}
P.X. Gallagher, {\em Some consequences of the Riemann hypothesis},
Acta Arith. \textbf{37} (1980), 339-343.

\bibitem{G}
S.M. Gonek, {\em On negative moments of the Riemann zeta-function},
Mathematika \textbf{36} (1989), 71-88.

\bibitem{G2}
S.M. Gonek, {\em The second moment of the reciprocal of the Riemann
zeta-function and its derivative}, Notes from a talk at MSRI,
summer 1999.

\bibitem{Ha}
C.B. Haselgrove, {\em A disproof of a conjecture of Polya},
Mathematika \textbf{5} (1958), 141-145.

\bibitem{HB}
D.R. Heath Brown, {\em The distribution and moments of the error term in the
Dirichlet divisor problem}, Acta Arith., \textbf{60} (4) (1992), 389-415.

\bibitem{Hej}
D. Hejhal, {\em On the distribution of log$|\zeta^{'}(\frac{1}{2} + it)|$
}, in {\em Number Theory, Trace Formula and Discrete Groups},
 eds. Aubert, Bombieri, and Goldfeld,
Academic Press, San Diego, 1989, 343-370.

\bibitem{Hl}
Edmund Hlawka, {\em The theory of uniform distribution}, AB
Academic Publishers, Berkhamsted, 1984.

\bibitem{HKO}
C.P. Hughes, J.P. Keating, and Neil O'Connell,
{\em Random matrix theory and the derivative of the Riemann zeta
function}, Proceedings of the Royal Society: A \textbf{456} (2000),
2611-2627.

\bibitem{In0}
A.E. Ingham, {\em The distribution of prime numbers},
Stechert-Hafner Service Agency, New York, 1964.

\bibitem{In}
A.E. Ingham, {\em On two conjectures in the theory of numbers},
Amer. J. Math., \textbf{64} (1942), 313-319.

\bibitem{Ivic}
Aleksandar Ivi\'{c}, {\em The theory of the Riemann zeta function
with applications}, John Wiley and Sons, New York, 1985.

\bibitem{KS}
J.P. Keating and N.C. Snaith,
{\em Random matrix theory and $\zeta(\frac{1}{2} + it)$},
Commun. Math. Phys., \textbf{214} (2000), 57-89.

\bibitem{Mo}
H.L. Montgomery, {\em The zeta function and prime numbers},
Proceedings of the Queen's Number Theory Conference, 1979,
Queen's Univ., Kingston, Ont., 1980, 1-31.

\bibitem{Ng}
N. Ng, {\em Limiting distributions and zeros of Artin
$L$-functions}, Ph.D. Thesis, University of British Columbia, fall
2000.

\bibitem{OR}
A.M. Odlyzko and H.J.J. te Riele, {\em Disproof of the Mertens
conjecture}, J. Reine Angew. Math., \textbf{357} (1985), 138-160.

\bibitem{Pr}
K. Prachar, {\em Primzahlverteilung}, Springer-Verlag, Berlin, 1957.

\bibitem{RS}
Michael Rubinstein and Peter Sarnak, {\em Chebyshev's Bias}, Journal
of Experiment. Math. \textbf{3} (1994), 173-197.

\bibitem{T}
E.C. Titchmarsh, {\em The theory of the Riemann zeta function,
Second Edition}, Oxford University Press, New York, 1986.

\end{thebibliography}
\end{document}